\documentclass[number, 1p]{elsarticle}
\usepackage{amsmath, amssymb}
\usepackage[all]{xy}
\usepackage{verbatim}
\usepackage{multirow}

\newcommand{\R}{\mathbb{R}}
\newcommand{\C}{\mathbb{C}}

\newtheorem{thm}{Theorem}[section]
\newtheorem{lem}[thm]{Lemma}
\newtheorem{prop}[thm]{Proposition}
\newtheorem{defn}[thm]{Definition}
\newtheorem{cor}[thm]{Corollary}
\newproof{pf}{Proof}

\begin{document} 
\author{Christopher Jankowski}
\address{Department of Mathematics, University of Pennsylvania, Philadelphia, PA 19104 USA}
\ead{cjankows@math.bgu.ac.il}
\fntext[fn1]{Present Address: Department of Mathematics,
Ben-Gurion University of the Negev, PO Box 653, Be'er Sheva 84105, Israel
\\ Fax: +972 8 647 7648}

\title{On type II$_0$ $E_0$-semigroups induced by boundary weight doubles}
\begin{abstract}
Powers has shown that each spatial $E_0$-semigroup can be obtained
from the boundary weight map of a $CP$-flow acting on $B(K \otimes
L^2(0, \infty))$ for some separable Hilbert space $K$. In this
paper, we define boundary weight maps through boundary weight
doubles $(\phi, \nu)$, where $\phi: M_n(\C) \rightarrow M_n(\C)$ is
a $q$-positive map and $\nu$ is a boundary weight over $L^2(0,
\infty)$. These doubles induce $CP$-flows over $K$ for $1<dim(K)<
\infty$ which then minimally dilate to $E_0$-semigroups by a theorem
of Bhat.  Through this construction, we obtain uncountably many
mutually non-cocycle conjugate $E_0$-semigroups for each $n>1, n \in
\mathbb{N}$.
\end{abstract}
\maketitle
\section{Introduction}
Let $H$ be a separable Hilbert space, denoting its inner product by
the symbol $( \ , \ )$ which is conjugate-linear in its first
coordinate and linear in its second.  A result of Wigner in
\cite{wigner} shows that every weakly continuous one-parameter group
of $*$-automorphisms $\{\alpha_t\}_{t \in \R}$ of $B(H)$ is
implemented by a strongly continuous unitary group $\{U_t\}_{t \in
\R}$ in that $\alpha_t(A)=U_tAU_t^*$ for all $A \in B(H)$ and $t \in
\R$. This leads us to pursue the more general task of classifying
all suitable semigroups of $*$-endomorphisms of $B(H)$:
\begin{defn}
We say a family $\{\alpha_t\}_{t \geq 0}$ of $*$-endomorphisms of
$B(H)$ is an \emph{$E_0$-semigroup} if:
\begin{enumerate}
\item $\alpha_{s+t}=\alpha_s \circ \alpha_t$ for all $s, t \geq 0$,
and $\alpha_0 (A)=A$ for all $A \in B(H)$.
\item  For each $f, g \in H$ and $A \in B(H)$, the inner product
$( f, \alpha_t (A) g)$ is continuous in $t$.
\item $\alpha_t (I)=I$ for all $t \geq 0$ (in other words, $\alpha$ is
\emph{unital}).
\end{enumerate}
\end{defn}

We have two different notions of what it means for two
$E_0$-semigroups to be the same, namely conjugacy and cocycle
conjugacy, the latter of which arises from Alain Connes' definition
of outer conjugacy.

\begin{defn}  Let $\alpha$ and $\beta$ be $E_0$-semigroups on $B(H_1)$ and
$B(H_2)$, respectively.  We say that $\alpha$ and $\beta$ are
\emph{conjugate} if there is a $*$-isomorphism $\theta$ from
$B(H_1)$ onto $B(H_2)$ such that $\theta \circ \alpha_t = \beta_t
\circ \theta$ for all $t \geq 0$.  We say that $\alpha$ and $\beta$
are \emph{cocycle conjugate} if $\alpha$ is conjugate to $\beta'$,
where $\beta'$ is an $E_0$-semigroup on $B(H_2)$ satisfying the
following condition: For some strongly continuous family of
unitaries $U=\{U_t: t \geq 0\}$ acting on $H_2$ and satisfying
$U_{t+s}=U_t \beta_t(U_s)$ for all $s, t \geq 0$, we have $\beta'_t
(A) = U_t \beta_t (A) U_t^*$ for all $A \in B(H_2)$ and $t \geq 0$.
Such a family of unitaries is called a \emph{unitary cocycle} for
$\beta$.
\end{defn}

$E_0$-semigroups are divided into three types based upon the
existence, and structure of, their units.  More specifically, let
$\alpha$ be an $E_0$-semigroup on $B(H)$.  A \textit{unit} for
$\alpha$ is a strongly continuous semigroup of bounded operators
$U=\{U(t): t \geq 0\}$ such that $\alpha_t(A)U(t) = U(t)A$ for all
$A \in B(H)$. Let $\mathcal{U}_\alpha$ be the set of all units for
$\alpha$. We say $\alpha$ is \textit{spatial} if $\mathcal{U}_\alpha
\neq \emptyset$, while we say that $\alpha$ is \textit{completely
spatial} if, for each $t \geq 0$, the closed linear span of the set
$\{ U_{1}(t_1) \cdots U_n(t_n) f: f \in H, t_i \geq 0 \textrm{ and }
U_i \in \mathcal{U}_\alpha \ \forall \ i, \sum t_i = t \}$ is $H$.
If an $E_0$-semigroup $\alpha$ is completely spatial, we say
it is of type I.  If $\alpha$ is spatial but is not completely
spatial, we say $\alpha$ is of type II.  If $\alpha$ has no units, we say
it is of type III. 

If $\alpha$ is of type I or II, we may further assign an integer $n
\in \mathbb{Z}_{\geq 0} \cup \{ \infty \}$ to $\alpha$, in which
case we say $\alpha$ is of type I$_n$ or II$_n$.  We call $n$ the
index of $\alpha$.  It was initially defined in different ways in
\cite{powersdiff} and \cite{arvindex}, and the connection between
these definitions was explored in \cite{powersprice}.  The index of
$\alpha$ is the dimension of a particular Hilbert space associated
to its units, and it is perhaps the most fundamental cocycle
conjugacy invariant for spatial $E_0$-semigroups.  Arveson showed in
\cite{arvindex} that the type I $E_0$-semigroups are entirely
classified (up to cocycle conjugacy) by their index: the type I$_0$
$E_0$-semigroups are semigroups of $*$-automorphisms, while for $n
\in \mathbb{N} \cup \{\infty\}$, every type I$_n$ $E_0$-semigroup is
cocycle conjugate to the CAR flow of rank $n$.

However, at the present time, we do not have such a classification
for those of type II or III.  The first type II and type III
examples were constructed by Powers in \cite{typeII} and
\cite{typeIII}.  Through Arveson's theory of
product systems, Tsirelson became the first to exhibit
uncountably many mutually non-cocycle conjugate $E_0$-semigroups
of types II and III (see \cite{T2}).  A dilation theorem of Bhat in \cite{Bhat} shows
that every unital $CP$-flow $\alpha$ can be dilated to an
$E_0$-semigroup, and that there is a minimal dilation $\alpha^d$ of
$\alpha$ which is unique up to conjugacy. Using Bhat's result,
Powers proved in \cite{hugepaper} that every spatial $E_0$-semigroup
can be obtained from the boundary weight map of a $CP$-flow over a
separable Hilbert space $K$.  In \cite{bigpaper}, he constructed
spatial $E_0$-semigroups using boundary weights over $K$ when
$dim(K)=1$ and then began to investigate the case when $dim(K)=2$.

Our goal is to use boundary weight maps to induce unital $CP$-flows
over $K$ for $1<dim(K)<\infty$ and to classify their minimal
dilations to $E_0$-semigroups up to cocycle conjugacy. To do so, we
define a natural boundary weight map $\rho \rightarrow \omega
(\rho)$ using a unital completely positive map $\phi$ and a
normalized boundary weight $\nu$ over $L^2(0, \infty)$.  The
necessary and sufficient condition that this map induce a unital
$CP$-flow $\alpha$ is that $\phi$ satisfies a definition of
$q$-positive analogous to that from \cite{hugepaper} (see Definition
\ref{qpos} and Proposition \ref{bdryweight}), in which case we say
that $\alpha$ is the $CP$-flow induced by the boundary weight double
$(\phi, \nu)$.  We develop a comparison theory for boundary weight
doubles $(\phi, \nu)$ and $(\psi, \nu)$ ($\phi$ and $\psi$ unital)
in the case that $\nu$ is a normalized unbounded boundary weight
over $L^2(0, \infty)$ of the form $\nu(\sqrt{I - \Lambda(1)} B
\sqrt{I - \Lambda(1)}) = (f,Bf)$, finding that the doubles induce
cocycle conjugate $E_0$-semigroups if and only if there is a
hyper maximal $q$-corner from $\phi$ to $\psi$ (see Definition
\ref{hyp} and Proposition \ref{hypqc}).

The problem of determining hyper maximal $q$-corners from $\phi$ to
$\psi$ becomes much easier if we focus on a particular class of
$q$-positive maps, called the $q$-pure maps, which have the least
possible $q$-subordinates (Definition \ref{qpure}). Given a
$q$-positive map $\phi$ acting on $M_n(\C)$ and a unitary $U \in
M_n(\C)$, we can form a new map $\phi_U$ by
$\phi_U(A)=U^*\phi(UAU^*)U$.  We describe the order isomorphism
between the $q$-subordinates of $\phi$ and those of $\phi_U$, which
in turn leads to the existence of a hyper maximal $q$-corner from
$\phi$ to $\phi_U$ if $\phi$ is unital and $q$-pure (Proposition
\ref{basischange}).   With this result in mind, we begin the task of
classifying the unital $q$-pure maps.  We find that the rank one
unital $q$-pure maps $\phi: M_n(\C) \rightarrow M_n(\C)$ are
precisely the maps $\phi(A)=\rho(A)I$ for faithful states $\rho$ on
$M_n(\C)$ (Proposition \ref{littlestates}). That these maps give us
an enormous class of mutually non-cocycle conjugate $E_0$-semigroups
in one of our main results (Theorem \ref{statesbig}).  Furthermore,
for $n>1$, none of the $E_0$-semigroups constructed from boundary
weight doubles satisfying the conditions of Theorem \ref{statesbig} are cocycle
conjugate to any of the $E_0$-semigroups obtained from
one-dimensional boundary weights by Powers in \cite{bigpaper} (Corollary
\ref{thesearenew}).

We turn our attention to the unital $q$-pure maps that are
invertible. These maps are best understood through their
(conditionally negative) inverses.  In Theorem \ref{phiu}, we find a
necessary and sufficient condition for an invertible unital map
$\phi$ on $M_n(\C)$ to be $q$-pure.  In this case, however, if $\nu$
is a normalized unbounded boundary weight of the form $\nu(\sqrt{I -
\Lambda(1)} B \sqrt{I - \Lambda(1)}) = (f,Bf)$, then the
$E_0$-semigroup induced by the boundary weight double $(\phi, \nu)$
is entirely determined by $\nu$.  This $E_0$-semigroup is the one
induced by $\nu$ in the sense of \cite{bigpaper}.

\section{Background}
\subsection{Completely positive maps}
Let $\phi: \mathfrak{U} \rightarrow \mathfrak{B}$ be a linear map
between $C^*$-algebras.  For each $n \in \mathbb{N}$, define
$\phi_n: M_n(\mathfrak{U}) \rightarrow M_n(\mathfrak{B})$ by

\begin{displaymath} \phi_n \left(\begin{array}{ccc} A_{11} & \cdots & A_{1n}
\\ \vdots & \ddots & \vdots \\ A_{n1} & \cdots & A_{nn}  \end{array}\right)
= \left(\begin{array}{ccc} \phi(A_{11}) & \cdots & \phi(A_{1n})
\\ \vdots & \ddots & \vdots \\ \phi(A_{n1}) & \cdots & \phi(A_{nn})
\end{array} \right). \end{displaymath}
We say that $\phi$ is completely positive if $\phi_n$ is positive
for all $n \in \mathbb{N}$.  A linear map $\phi: B(H_1)
\rightarrow B(H_2)$ is completely positive if and only if for all
$A_1, \ldots A_n \in B(H_1)$, $f_1, \ldots, f_n \in H_2$, and $n \in
\mathbb{N}$, we have
$$\sum_{i,j=1}^n (f_i, \phi(A_i^*A_j)f_j) \geq 0.$$ 
Stinespring's Theorem asserts that if
$\mathfrak{U}$ is a unital $C^*$-algebra and $\phi: \mathfrak{U}
\rightarrow B(H)$ is a unital completely positive map, then $\phi$
dilates to a $*$-homomorphism in that there is a Hilbert space $K$,
a $*$-homomorphism $\pi: \mathfrak{U} \rightarrow B(K)$, and an
isometry $V: H \rightarrow K$ such that
$$\phi(A)=V^*\pi(A)V$$ for all $A \in \mathfrak{U}$. 

From the work of Choi (\cite{choi}) and Arveson (\cite{arveson}), 
we know that a normal linear map $\phi: B(H_1)
\rightarrow B(H_2)$ is completely positive if and only if it can be
written in the form
$$\phi(A) = \sum_{i=1}^n S_i A S_i^*$$ for 
some $n \in \mathbb{N} \cup \{\infty\}$
and maps $S_i: H_1
\rightarrow H_2$ which are linearly independent over
$\ell_2(\mathbb{N})$ in the sense that if $\sum_{i=1}^{r \leq n} z_i S_i = 0$ 
for a sequence $\{z_i\}_{i=1}^r \in
\ell_2(\mathbb{N})$, then $z_i=0$ for all $i$.  With these
hypotheses satisfied, the number $n$ is unique.  We will use the
above conditions for complete positivity interchangeably.

\subsection{Conditionally negative maps}
We say a self-adjoint linear map $\psi: B(K) \rightarrow B(K)$ is
\emph{conditionally negative} if, whenever $\sum_{i=1}^m
A_if_i=0$ for $A_1, \ldots, A_m \in B(K)$, $f_1, \ldots, f_m \in K$,
and $m \in \mathbb{N}$, we have $\sum_{i=1}^m (f_i,
\psi(A_i^*A_j)f_j) \leq 0$.  If $K=\C^n$, then from the literature
(see, for example, Theorem 3.1 of \cite{pownxn}) we
know that $\psi$ has the form
$$\psi(A) = sA + YA + AY^* - \sum_{i=1}^p \lambda_i S_iAS_i^*,$$ where $s
\in \R$, $tr(Y)=0$, and for all $i$ and $j$ we have $\lambda_i
>0$, $tr(S_i)=0$ and $tr(S_i^*S_j)=n \delta_{ij}$, where
$p \leq n^2$ is independent of the maps $S_i$.

This form for $\psi$ is unique in the sense that if $\psi$ is
written in the form
$$\psi(A)= tA + ZA +AZ^* - \sum_{i=1}^p \mu_i T_i A T_i^*,$$ where $t \in
\R$, $tr(Z)=0$, and for all $i$ and $j$ we have $\mu_i >0$,
$tr(T_i)=0$, and $tr(T_i^*T_j)= n \delta_{ij}$, then $s=t$, $Z=Y$,
and $\sum_{i=1}^p \lambda_i S_iAS_i^*= \sum_{i=1}^p \mu_i T_i A
T_i^*$ for all $A \in M_n(\C)$. Indeed, let $\{v_k\}_{k=1}^n$ be any
orthonormal basis for $\C^n$, let $h_k = v_k / \sqrt{n}$ for each
$k$, let $f \in \C^n$ be arbitrary, and for $k=1, \ldots, n,$ define
$A_k \in M_n(\C)$ by $A_k = f h_k^*$. Using the trace conditions, we
find
\begin{eqnarray*} \sum_{k=1}^n \psi(A_k)h_k & = & \sum_{k=1}^n (h_k, h_k)sf
+ \sum_{k=1}^n (h_k, h_k)Yf +   \sum_{k=1}^n (h_k, Y^* h_k)f \\
& \ & - \sum_{k=1}^n \Big(\sum_{i=1}^p \lambda_i(h_k, S_i^*h_k)S_i
f\Big)
\\
&= & sf + Yf + 0 - \sum_{i=1}^p \Big( \sum_{k=1}^n \lambda_i(h_k,
S_i^* h_k)
S_i f \Big) \\
& = & sf + Yf - \sum_{i=1}^p \lambda_i (0) S_i f = sf + Yf.
\end{eqnarray*}
An analogous computation shows that $\sum_{k=1}^n \psi(A_k)h_k = tf
+ Zf$. Since $f \in \C^n$ was arbitrary, we conclude $(t-s)I= Y-Z$.
Therefore, $tr((t-s)I) = tr(Y-Z)=0$, so $t=s$ and $Y=Z$.
Consequently, $\sum_{i=1}^p \lambda_i S_iAS_i^*= \sum_{i=1}^p \mu_i
T_i A T_i^*$ for all $A \in M_n(\C)$.

\subsection{$CP$-flows and Bhat's theorem}  Let $K$ be a separable Hilbert
space and let $H= K \otimes L^2(0, \infty)$.  We identify $H$ with
$L^2((0, \infty); K)$, the space of $K$-valued measurable functions
on $(0, \infty)$ which are square integrable. Under this
identification, the inner product on $H$ is
$$(f,g) = \int_0 ^\infty (f(x), g(x)) dx.$$  Let
$U=\{U_t\}_{t \geq 0}$ be the right shift semigroup on $H$, so for 
all $t \geq 0$ and $f \in H$ we have $(U_t
f)(x) = f(x-t)$ for $x>t$ and $(U_t f)(x)=0$ otherwise.  Let
$\Lambda: B(K) \rightarrow B(H)$ be the map defined by
$(\Lambda(A)f)(x)= e^{-x}Af(x)$ for all $A \in B(K), f \in H$.  

\begin{defn}  Assume the above notation.
A strongly continuous semigroup $\alpha=\{\alpha_t: t \geq 0\}$ of
completely positive contractions of $B(H)$ into itself is a
\textit{CP-flow} if $\alpha_t(A)U_t = U_tA$ for all $A \in B(H)$.
\end{defn}

A theorem of Bhat in \cite{Bhat} allows us to generate
$E_0$-semigroups from unital $CP$-flows, and, more generally, from
strongly continuous completely positive semigroups of unital maps on
$B(H)$, called $CP$-semigroups. We give a reformulation of Bhat's
theorem (see Theorem 2.1 of \cite{bigpaper}):

\begin{thm}\label{dilation} Suppose $\alpha$ is a unital $CP$-semigroup of $B(H_1)$.
Then there is an $E_0$-semigroup $\alpha^d$ of $B(H_2)$ and an
isometry $W: H_1 \rightarrow H_2$ such that
$$\alpha_t(A)=W^*\alpha_t ^d(WAW^*)W$$ and $\alpha_t(WW^*) \geq
WW^*$ for all $t > 0$.  If the projection $E=WW^*$ is minimal in
that the closed linear span of the vectors
$$\alpha_{t_1}^d(EA_1E) \cdots \alpha_{t_n}^d(EA_nE)Ef$$ for $f \in
K, A_i \in B(H_1)$ and $t_i \geq 0$ for all $i=1, 2, \ldots, n$ and
$n=1, 2, \ldots$ is $H_2$,
then $\alpha^d$ is unique up to conjugacy.
\end{thm}

In \cite{hugepaper}, Powers showed that every spatial
$E_0$-semigroup acting on $B(\mathfrak{H})$ (for $\mathfrak{H}$ a
separable Hilbert space) is cocycle conjugate to an $E_0$-semigroup
which is a $CP$-flow, and that every $CP$-flow over $K$ arises from
a \textit{boundary weight map} over $H=K \otimes L^2(0, \infty)$.
The boundary weight map $\rho \rightarrow \omega (\rho)$ of a
$CP$-flow $\alpha$ associates to every $\rho \in B(K)_*$ a boundary
weight, that is, a linear functional $\omega(\rho)$ acting on the
null boundary algebra
$$\mathfrak{A}(H)= \sqrt{I_H - \Lambda(I_K)}B(H)\sqrt{I_H-\Lambda(I_K)}
$$ which is normal in the following sense:  If we define a linear
functional $\ell(\rho)$ on $B(H)$ by
$$\ell(\rho)(A)=\omega(\rho)\Big(\sqrt{I_H-\Lambda(I_K)}A \sqrt{I_H -
\Lambda(I_K)}\Big),$$ then $\ell(\rho) \in B(H)_*$.  If
$\omega(\rho)(I_H-\Lambda(I_K)) = \rho(I_K)$ for all $\rho \in
B(K)_*$, then $\alpha$ is unital.  For the sake of neatness, we will
omit the subscripts $H$ and $K$ from the previous sentence when
they are clear.  Let $\delta$ be the generator of $\alpha$, and define $\Gamma:
B(H) \rightarrow B(H)$ by $\Gamma(A)= \int_0 ^ \infty e^{-t}U_t AU_t^*$.
The resolvent $R_\alpha:= (I - \delta)^{-1}$ of $\alpha$ satisfies $R_\alpha(A) = \int_0 ^\infty e^{-t}\alpha_t(A) dt$
for all $A \in B(H)$.
Its associated predual map $\hat{R}_\alpha$ is given by
\begin{equation}\label{resolvent}
\hat{R}_\alpha(\eta) = \hat{\Gamma}(\omega(\hat{\Lambda}\eta) + \eta)
\end{equation}
for all $\eta \in B(H)_*$.

A $CP$-flow $\alpha$ over $K$ is entirely determined by a set of normal
completely positive contractions $\pi^{\#}=\{\pi^{\#}_t: t > 0\}$ from
$B(H)$ into $B(K)$,
called the \textit{generalized boundary representation} of $\alpha$.
Its relationship to the boundary weight map is as follows.  For
each $t>0$, denote by $\hat{\pi}_t: B(K)_* \rightarrow B(H)_*$ 
the predual map induced by $\pi_t^\#$.  For the
truncated boundary weight maps $\rho \rightarrow \omega_t(\rho) \in
B(H)_*$ defined by
\begin{eqnarray}\label{truncate} \omega_t(\rho)(A)=
\omega(\rho)\Big(U_tU_t^*AU_tU_t^*\Big),\end{eqnarray} we have
$\hat{\pi}_t=\omega_t(I + \hat{\Lambda}\omega_t)^{-1}$ and $\omega_t
= \hat{\pi}_t(I - \hat{\Lambda} \hat{\pi}_t)^{-1}$ for all $t>0$.
The maps $\{\pi_b^\#\}_{b>0}$ have a $\sigma$-strong limit
$\pi_0^\#$ as $b \rightarrow 0$ for each $A \in \bigcup_{t>0} U_t
B(H) U_t^*$, called the \textit{normal spine} of $\alpha$. If
$\alpha$ is unital, then the index of $\alpha^d$ as an
$E_0$-semigroup is equal to the rank of $\pi_0^\#$ as a completely
positive map (Theorem 4.49 of \cite{hugepaper}).  

Having seen that every $CP$-flow has an associated boundary weight map, we
would like to approach the situation from the opposite direction.
More specifically, under what conditions is a map 
$\rho \rightarrow \omega(\rho)$ from $B(K)_*$ to weights acting on $\mathfrak{A}(H)$ the boundary
weight map of a $CP$-flow over $K$?  Powers has found the answer (see Theorem 3.3 
of \cite{bigpaper}):
\begin{thm}\label{powersthm}
If $\rho \rightarrow \omega(\rho)$ is a completely positive mapping
from $B(K)_*$ into weights on $B(H)$ satisfying
$\omega(\rho)(I-\Lambda(I_K)) \leq \rho(I_K)$ for all positive
$\rho \in B(K)_*$, and if the maps $\hat{\pi}_t : = \omega_t(I +
\hat{\Lambda}\omega_t)^{-1}$ are completely positive contractions
from $B(K)_*$ into $B(H)_*$ for all $t >0$, then $\rho \rightarrow
\omega(\rho)$ is the boundary weight map of a $CP$-flow over $K$.
The $CP$-flow is unital if and only if $\omega(\rho)(I-\Lambda(I_K))
= \rho(I_K)$ for all $\rho \in B(K)_*$.
\end{thm}

If $dim(K)=1$, the boundary weight map is just $c \in \C \rightarrow
\omega(c) = c \omega(1)$, so we may view our boundary weight map as
a single positive boundary weight $\omega:= \omega(1)$ acting on
$\mathfrak{A}(L^2(0, \infty))$.  Since the functional $\ell$ defined
on $B(H)$ by $$\ell(A)=\omega\Big(\sqrt{I - \Lambda(1)} A \sqrt{I -
\Lambda(1)}\Big)$$ is positive and normal, it has the form $\ell(A)
= \sum_{k=1}^n (f_k, A f_k)$ for some mutually orthogonal vectors
$\{f_k\}_{k=1}^{n \in \mathbb{N} \cup \{\infty\}}$, so
$$\omega \Big(\sqrt{I - \Lambda(1)}A\sqrt{I - \Lambda(1)}\Big) =
\sum_{k=1}^n (f_k, A f_k)$$ for all $A \in B(H)$.  If $\omega$ is
\textit{normalized} (that is, $\omega(I-\Lambda(1))=1$), then
$\sum_{k=1}^n ||f_k||^2 = 1$.  In \cite{bigpaper}, Powers induced
$E_0$-semigroups using normalized boundary weights over $L^2(0,
\infty)$.   The \ type \ of \ $E_0$-semigroup \ $\alpha^d$ 
\ induced \ by \ a \
normalized \ boundary \ weight $\omega(\sqrt{I - \Lambda(1)}A \sqrt{I -
\Lambda(1)}) = \sum_{k=1}^n (f_k, A f_k)$ depends on whether
$\omega$ is bounded in the sense that for some $r>0$ we have
$|\omega(B)| \leq r ||B||$ for all $B \in \mathfrak{A}(H)$.  Results
from \cite{hugepaper} imply that $\alpha^d$ is of type I$_n$ if
$\omega$ is bounded and of type II$_0$ if $\omega$ is unbounded.  If
$\omega$ is unbounded, then both $\omega_t(I)$ and $\omega_t(\Lambda(1))$ approach infinity
as $t$ approaches zero.  We will focus on normalized unbounded boundary weights over $L^2(0,
\infty)$ of the form $\omega(\sqrt{I - \Lambda(1)}A \sqrt{I -
\Lambda(1)})=(f,Af).$ We note that, as discussed in detail in
\cite{markie}, such boundary weights are not normal weights.

If $\alpha$ and $\beta$ are $CP$-flows, we say that $\alpha \geq \beta$ if $\alpha_t 
- \beta_t$ is completely positive for all $t \geq 0$. The subordinates of a 
$CP$-flow are entirely determined by the
subordinates of its generalized boundary representation (see Theorem 3.4
of \cite{bigpaper}):

\begin{thm} \label{bdryrep} Let $\alpha$ and $\beta$ be $CP$-flows over $K$ with
generalized boundary representations $\pi^\#=\{\pi^{\#}_t\}$ and
$\xi^\# =\{\xi^{\#}_t\}$, respectively. Then $\beta$ is subordinate
to $\alpha$ if and only if
$\pi^{\#}_t - \xi^{\#}_t$ is completely positive for all $t>0$.
\end{thm}

Given 
two unital $CP$-flows $\alpha$ and $\beta$, it is natural to ask when their
minimally dilated $E_0$-semigroups are cocycle conjugate.  The following definition
from \cite{hugepaper} provides us with a key:
\begin{defn}
Let $\alpha$ and $\beta$ be $CP$-flows over $K_1$ and $K_2$, respectively, where
$H_1 = K_1 \otimes L^2(0, \infty)$ and $H_2= K_2 \otimes L^2(0, \infty)$.  We say
that a family of linear maps
$\gamma=\{\gamma_t: t \geq 0\}$ from $B(H_2, H_1)$ into itself is a
flow corner from $\alpha$ to $\beta$ if the family of maps
$\Theta=\{\Theta_t: t \geq 0\}$ defined by
\begin{displaymath} \Theta_t \left( \begin{array}{cc} A_{11} &  A_{12}  \\
A_{21} & A_{22}  \\
\end{array} \right)=
\left( \begin{array}{cc} \alpha_t(A_{11}) &  \gamma_t(A_{12})  \\
\gamma_t^*(A_{21}) & \beta_t(A_{22})  \\
\end{array} \right)
\end{displaymath} is a $CP$-flow over $K_1 \oplus K_2$.

If $\gamma$ is a flow corner from $\alpha$ to $\beta$, we consider
subordinates $\Theta'=\{\Theta_t': t \geq 0\}$ of $\Theta$ that are $CP$-flows of the form
\begin{displaymath} \Theta_t' \left( \begin{array}{cc} A_{11} &  A_{12}  \\
A_{21} & A_{22}  \\
\end{array} \right):=
\left( \begin{array}{cc} \alpha'_t(A_{11}) &  \gamma_t(A_{12})  \\
\gamma_t^*(A_{21}) & \beta'_t(A_{22})  \\
\end{array} \right).
\end{displaymath}  We say that $\gamma$ is a hyper
maximal flow corner from $\alpha$ to $\beta$ if, for every such
subordinate $\Theta'$ of $\Theta$, we have $\alpha=\alpha'$ and
$\beta=\beta'$.
\end{defn}

Our results will involve type II$_0$ $E_0$-semigroups.  These are spatial
$E_0$-semigroups which are not semigroups of $*$-automorphisms
and have only one unit $V=\{V_t\}_{t \geq 0}$ up to scaling by $e^{t \lambda}$ for
$\lambda \in \C$.
In the case that unital $CP$-flows $\alpha$ and $\beta$
minimally dilate to type II$_0$ $E_0$-semigroups, we have a necessary and
sufficient condition for $\alpha^d$ and $\beta^d$ to be cocycle
conjugate (Theorem 4.56 of \cite{hugepaper}):

\begin{thm}\label{hyperflowcorn}  Suppose $\alpha$ and $\beta$ are unital
$CP$-flows over
$K_1$ and $K_2$ and $\alpha^d$ and $\beta^d$ are their minimal dilations
to $E_0$-semigroups.  Suppose $\gamma$ is a hyper maximal flow corner from $\alpha$
to $\beta$.  Then $\alpha^d$ and $\beta^d$ are cocycle conjugate.  Conversely,
if $\alpha^d$ is a type II$_0$ and $\alpha^d$ and $\beta^d$ are cocycle conjugate,
then there is a hyper maximal flow corner from $\alpha$ to $\beta$.

\end{thm}

We will later use this theorem
to determine a necessary and sufficient condition for some of the
$E_0$-semigroups we construct to be cocycle conjugate (see Definition \ref{hyp} and
Proposition \ref{hypqc}).

\section{Our boundary weight map}
Recall that a completely positive linear map $\phi$ can have
negative eigenvalues.  Moreover, even if $I + t \phi$ is invertible
for a given $t$, it does not necessarily follow that $\phi(I + t
\phi)^{-1}$ is completely positive.  In our boundary weight
construction, we will require a special kind of completely positive
map:

\begin{defn}\label{qpos}  A linear map $\phi: M_n(\C) \rightarrow
M_n(\C)$ is
\emph{$q$-positive} if $\phi$ has no negative eigenvalues and
$\phi(I + t \phi)^{-1}$ is completely
positive for all $t \geq 0$.
\end{defn}

Henceforth, we naturally identify a finite-dimensional Hilbert space
$K$ with $\C^n$ and $B(K \otimes L^2(0, \infty))$ with $M_n(B(L^2(0,
\infty)))$.  Under these identifications, the right shift $t$ units
on $K \otimes L^2(0, \infty)$ is the matrix whose $ij$th entry is
$\delta_{ij} V_t$ for $V_t$ the right shift on $L^2(0, \infty)$. The
map $\Lambda_{n \times n}: B(K) \rightarrow B(K \otimes L^2(0,
\infty))$ sends an $n \times n$ matrix $B=(b_{ij}) \in M_n(\C)$ to
the matrix $\Lambda_{n \times n}(B)$ whose $ij$th entry
is $b_{ij} \Lambda(1) \in B(L^2(0, \infty))$.  The null boundary
algebra $\mathfrak{A}(H)$ is simply $M_n (\mathfrak{A}(L^2(0,
\infty)))$.

Given a boundary weight $\nu$ over $L^2(0, \infty)$, we write
$\Omega_{\nu, n \times k}$ for the map that sends an $n \times k$
matrix $A=(A_{ij}) \in M_{n \times k}(\mathfrak{A}(L^2(0, \infty)))$
to the matrix $\Omega_{\nu, n \times k} (A) \in M_{n \times k}(\C)$
whose $ij$th entry is $\nu(A_{ij})$.  We will suppress the
integers $n$ and $k$ when they are clear, writing the above maps as
$\Omega_\nu$ and $\Lambda$. In the proposition and corollary that
follow, we show how to construct a $CP$-flow using a $q$-positive
map $\phi: M_n(\C) \rightarrow M_n(\C)$, a normalized boundary
weight $\nu$ over $L^2(0, \infty)$, and the map $\Omega_\nu:=
\Omega_{\nu, n \times n}: \mathfrak{A}(H) \rightarrow M_n(\C)$. The
map $\Omega_\nu$ is completely positive since $\nu$ is positive.

\begin{prop} \label{bdryweight} Let $H =\C^n \otimes L^2(0, \infty)$.
Let $\phi: M_n(\C) \rightarrow M_n(\C)$ be a unital completely
positive map with no negative eigenvalues, and let $\nu$ be a normalized unbounded boundary
weight over $L^2(0, \infty)$.  Then the map $\rho \rightarrow
\omega(\rho)$ from $M_n(\C)^*$ into boundary weights on
$\mathfrak{A}(H)$ defined by
$$\omega (\rho) (A) = \rho(\phi(\Omega_\nu(A))).$$
is completely positive.  Furthermore, the maps $\hat{\pi}_t:= \omega_t (I + \hat{\Lambda}\omega_t)^{-1}$
define normal completely positive contractions $\pi_t ^\#$ of $B(H)$ into $M_n(\C)$
for all $t
> 0$ if and only if $\phi$ is $q$-positive.
\end{prop}
\begin{pf}  
The map $\rho \rightarrow \omega(\rho)$ is completely positive since
it is the composition of two completely positive maps. Before
proving either direction, we let $s_t= \nu_t(\Lambda(1))$
for all $t>0$ and prove the equality
\begin{eqnarray} \label{pi} \hat{\pi}_t(\rho) = \rho \Big(\phi(I + s_t \phi)^{-1}
\Omega_{\nu_t}\Big) \end{eqnarray} for all $\rho \in M_n(\C)^*$. Denoting by $U_t$ 
the right shift on $H$ for every $t>0$, we claim
that $(I+\hat{\Lambda}\omega_t)^{-1} = (I + s_t \hat{\phi})^{-1}$.
Indeed, for arbitrary $t>0$, $B \in M_n(\C)$, and $\rho \in
M_n(\C)^*$, we have
$$\hat{\Lambda}\omega_t(\rho)(B)=
\rho\Big(\phi\Big(\Omega_\nu(U_tU_t^*\Lambda(B)U_tU_t^*)\Big)\Big)=
\rho\Big(\phi\Big(\Omega_{\nu_t}\Big(\Lambda(B)\Big)\Big)\Big)=s_t \rho(\phi(B)),$$
hence $\hat{\Lambda}\omega_t = s_t \hat{\phi}$ and
$(I+\hat{\Lambda}\omega_t)^{-1} = (I + s_t \hat{\phi})^{-1}$. 

For any $t>0$ and $A \in B(H)$, we have
\begin{eqnarray*} \hat{\pi}_t(\rho)(A) & = &\omega_t
(I + \hat{\Lambda}\omega_t)^{-1}(\rho)(A) = \Big((I +
\hat{\Lambda}\omega_t)^{-1}(\rho)\Big)(\phi(\Omega_{\nu_t}(A)))
\\ & = & \Big((I + s_t \hat{\phi})^{-1}(\rho)\Big)(\phi(\Omega_{\nu_t}(A)))
= \rho\Big((I + s_t \phi)^{-1} \phi(\Omega_{\nu_t}(A))\Big)\\ & = &
\rho\Big(\phi(I + s_t \phi)^{-1}(\Omega_{\nu_t}(A))\Big),
\end{eqnarray*} establishing \eqref{pi}.

Assume the hypotheses of the backward direction and let $t>0$.  By construction, 
$\hat{\pi}_t$ maps $M_n(\C)^*$ into $B(H)_*$.  It is also a contraction,
since for all $\rho \in M_n(\C)^*$ we have
\begin{eqnarray*} ||\hat{\pi}_t(\rho)|| & = & \Big| \Big|\rho\Big(\phi(I +
s_t\phi)^{-1}\Omega_{\nu_t}\Big)\Big| \Big|  \leq  ||\rho|| \ ||\phi(I + s_t \phi)^{-1}
\Omega_{\nu_t}|| \\ & = & ||\rho|| \ ||\phi(I + s_t \phi)^{-1}
\Omega_{\nu_t}(I)|| 
= ||\rho|| \ \Big|\Big| \phi(I + s_t \phi)^{-1}\Big(\nu_t(I)I_{\C^n}\Big)
\Big|\Big| \\ & = & ||\rho|| \Big| \Big|\frac{\nu_t(I)}{1+s_t} I_{\C^n} \Big|\Big| =
||\rho || \frac{\nu_t(I)}{1+\nu_t(\Lambda(1))}\leq ||\rho||,
\end{eqnarray*}
where the last inequality follows from the fact that
$$\nu_t(I-\Lambda(1)) \leq
\nu(I-\Lambda(1)) = 1.$$  Therefore, for every $t>0$, $\hat{\pi}_t$ defines
a normal contraction $\pi_t ^\#$ from $B(H)$ into $M_n(\C)$
satisfying $\hat{\pi}_t(\rho) = \rho \circ \pi_t ^\#$ for all $\rho \in M_n(\C)^*$. 
From Eq.\eqref{pi}
we see $\pi_t ^\# = \phi(I + s_t \phi)^{-1} \Omega_{\nu_t}$, so $\pi_t^\#$ is 
the composition of completely positive maps and is thus completely positive for all $t>0$.

Now assume the hypotheses of the forward direction.  By
unboundedness of $\nu$, the (monotonically decreasing) values
$\{s_t\}_{t>0}$ form a set equal to either $(0, \infty)$ or $[0,
\infty)$.  Choose any $t>0$ such that $s_t>0$. Let $T \in B(H)$ be
the matrix with $ij$th entry $(1/\nu_t(I))I$, and let
$\kappa_t: M_n(\C) \rightarrow B(H)$ be the map that sends
$B=(b_{ij}) \in M_n(\C)$ to the matrix $\kappa_t(B) \in B(H)$ whose
$ij$th entry is $(b_{ij}/\nu_t(I)) I$. We note that
$\kappa_t$ is the Schur product $B \rightarrow B \cdot T$, which is
completely positive since $T$ is positive. For all $B \in M_n(\C)$,
we have
$$\phi(I + s_t \phi)^{-1}(B)= \pi^\#_t (\kappa_t(B)),$$
so $\phi(I + s_t \phi)^{-1}$ is the composition of completely
positive maps and is thus completely positive.  As noted above, the
values $\{s_t\}_{t>0}$ span $(0, \infty)$, so $\phi$ is
$q$-positive. \qed \end{pf}

\begin{cor}\label{bdryweightcor}  The map $\rho \rightarrow
\omega(\rho)$ in Proposition \ref{bdryweight} is the boundary weight
map of a unital $CP$-flow $\alpha$ over $\C^n$, and the Bhat minimal
dilation $\alpha^d$ of $\alpha$ is a type II$_0$ $E_0$-semigroup.
\end{cor}

\begin{pf} The first claim of the corollary follows immediately
from Theorem \ref{powersthm} and Proposition \ref{bdryweight} since
\begin{equation}\label{tech} \omega(\rho)(I-\Lambda(I_{\C^n}))=\rho(\phi(I_{\C^n}))
=\rho(I_{\C^n})
\end{equation} for all $\rho \in M_n(\C)^*$.  For the second
assertion, we note that by Theorem 4.49 of \cite{hugepaper}, the
index of $\alpha^d$ is equal to the rank of the normal spine
$\pi_0^{\#}$ of $\alpha$, where $\pi_0^\#$ is the $\sigma$-strong
limit of the maps $\{\pi_b^\#\}_{b>0}$ for each $A \in \bigcup_{t>0}
U_tB(H)U_t^*$. Fix $t>0$, and let $A \in U_tB(H)U_t^*$. From formula
\eqref{pi},
$$\pi_b^{\#}(A) = \phi(I + \nu_b(\Lambda(1)) \phi)^{-1}
(\Omega_{\nu_b}(A)).$$  For all $b<t$ we have
$||\Omega_{\nu_b}(A)||=||\Omega_{\nu_t}(A)||< \infty$.  Since
$\nu_b(\Lambda(1)) \rightarrow \infty$ as $b \rightarrow 0$, we
conclude $\lim_{b \rightarrow 0} ||\pi_b^\# (A)|| =0,$  hence
$\pi_0^{\#} = 0$ and the index of $\alpha$ is zero.  However,
$\alpha^d$ is not completely spatial since $\alpha$ is 
not derived from the zero boundary weight map (see Lemma 4.37 and Theorem
4.52 of \cite{hugepaper}), so $\alpha^d$ is of type II$_0$. \qed \end{pf}

Given a $q$-positive $\phi: M_n(\C) \rightarrow M_n(\C)$ and a
normalized unbounded boundary weight $\nu$ over $L^2(0, \infty)$, we
call $(\phi, \nu)$ a \textit{boundary weight double}. As we have
seen, if $\phi$ is unital then the boundary weight double naturally
defines a boundary weight map through the construction of
Proposition \ref{bdryweight}, inducing a type II$_0$ $E_0$-semigroup
$\alpha^d$ which is unique up to conjugacy by Theorem
\ref{dilation}. We should note that it is not necessary for $\phi$
to be unital in order for the boundary weight double to induce a
$CP$-flow: If $\phi$ is any $q$-positive contraction such that $||\nu_t(I) \phi(I +
\nu_t(\Lambda(1)) \phi)^{-1}|| \leq 1$ for all $t
> 0$, then the arguments given in the proofs of Proposition
\ref{bdryweight} and Corollary \ref{bdryweightcor} show that the
boundary weight double $(\phi, \nu)$ induces a $CP$-flow $\alpha$.
However, if $\phi$ is not unital, then by Eq.\eqref{tech} and
Theorem \ref{powersthm}, neither is $\alpha$.

Motivated by \cite{hugepaper}, we make the following definition:
\begin{defn}  Suppose $\alpha: B(H_1) \rightarrow B(K_1)$ and $\beta: B(H_2) \rightarrow
B(K_2)$ are normal and completely positive.  Write each $A \in B(H_1 \oplus H_2)
$ as $A=(A_{ij})$, where $A_{ij} \in B(H_j, H_i)$ for each
$i,j=1,2$.  We say a linear map $\gamma:
B(H_2, H_1) \rightarrow B(K_2, K_1)$ is a corner from $\alpha$ to
$\beta$ if $\psi: B(H_1 \oplus H_2) \rightarrow B(K_1 \oplus K_2)$
defined by
\begin{displaymath} \psi \left( \begin{array}{cc} A_{11} & A_{12} \\
A_{21} & A_{22} \end{array}  \right) = \left(
\begin{array}{cc} \alpha(A_{11}) & \gamma(A_{12}) \\ \gamma^*(A_{21})
& \beta(A_{22}) \end{array} \right)
\end{displaymath} is a normal completely positive map. \end{defn}

We will repeatedly use the following lemma, which gives us the form of any corner between
normal completely positive contractions of finite index.  We believe 
that this result is already present in the literature, but 
we present a proof here for the sake of completeness:

\begin{lem}\label{corners}  Let $H_1, H_2, K_1,$ and $K_2$
be separable Hilbert spaces.  Let $\alpha: B(H_1) \rightarrow
B(K_1)$ and $\beta: B(H_2) \rightarrow B(K_2)$ be normal completely
positive contractions of the form
$$\alpha(A_{11})= \sum_{i=1}^n S_i A_{11}S_i^*,  \beta(A_{22}) =
\sum_{j=1}^p T_j A_{22} T_j^*,$$ where $n,p \in \mathbb{N}$ and the
sets of maps $\{S_i\}_{i=1}^n$ and $\{T_j\}_{j=1}^p$ are both
linearly independent. A linear map $\gamma: B(H_2, H_1) \rightarrow
B(K_2,K_1)$ is a corner from $\alpha$ to $\beta$ if and only if for
all $A_{12} \in B(H_2,H_1)$ we have
$$\gamma(A_{12})= \sum_{i,j} c_{ij} S_i A_{12} T_j^*,$$ where
$C=(c_{ij}) \in M_{n \times p}(\C)$ is any matrix such that $||C||
\leq 1$.
\end{lem}

\begin{pf} For the backward direction, let $C=(c_{ij}) \in M_{n
\times p}(\C)$ be any contraction, and define a linear map $\gamma:
B(H_2, H_1) \rightarrow B(K_2,K_1)$ by $\gamma(A) = \sum_{i,j}
c_{ij} S_iAT_j^*$.  We need to show that the map
\begin{displaymath}L \left( \begin{array}{cc} A_{11} & A_{12}
\\ A_{21} & A_{22} \end{array} \right)=
\left( \begin{array}{cc} \alpha(A_{11}) &
\gamma(A_{12}) \\ \gamma^*(A_{21}) & \beta(A_{22}) \end{array} \right)
\end{displaymath} is normal and completely positive.  To prove this, we first assume
that $n \geq p$ and note that by Polar Decomposition we may write
$C_{n \times p} = V_{n \times p} T_{p \times p}$, where $V_{n \times
p}$ is a partial isometry of rank $p$ and $T$ is positive. Unitarily
diagonalizing $T$ we see $C_{n \times p} = V_{n \times p} W_{p
\times p}^* D_{p \times p} W_{p \times p}$. We may easily add
columns to $V_{n \times p} W_{p \times p}^*$ to form a unitary
matrix in $M_n(\C)$, which we call $U^*$. Defining $\tilde{D}
= (d_{ij}) \in M_{n \times p}(\C)$
to be the matrix obtained from $D$ by adding $n-p$ rows of zeroes,
we see $U^*\tilde{D} = V_{n \times p} W_{p \times p}^*D$, so $C_{n
\times p} = U^* \tilde{D} W_{p \times p}$ and
$$UC_{n \times p}W_{p \times p}^*=\tilde{D}.$$ In other words, \begin{displaymath} 
 \sum_{i,j}
c_{ij} u_{ki} \overline{w_{\ell j}} = \left\{
\begin{array}{ll} \delta_{k \ell}
d_{k \ell} & \textrm{ if } k \leq p \\
0 & \textrm{ if } k > p \\ \end{array} \right \}.
\end{displaymath}
Next, define $\{S_i'\}_{i=1}^n: H_1 \rightarrow K_1$ and
$\{T_j\}_{j=1}^p: H_2 \rightarrow K_2$ by
$$S_i ' = \sum_{k=1}^n \overline{u_{ik}} S_{k}, \ \ T_j ' =
\sum_{\ell =1}^p \overline{w_{j \ell}} T_j,$$ so $S_i = \sum_{k=1}^n
{u_{ki}} S_k'$ and $T_j = \sum_{\ell=1}^p w_{\ell j} T_j'$ for all
$i$ and $j$.  \\ \\ Since $U$ and $W$ are unitary, it follows that $||D||=||C|| \leq 1$ and that the maps
$\{S_i'\}_{i=1}^n$ are linearly independent, as are the maps
$\{T_j'\}_{j=1}^p$.
We observe that for any $A_{11} \in B(H_1)$ and $A_{22} \in B(H_2)$,
$$\sum_{i=1}^n S_i A_{11} S_i ^* = \sum_{i=1}^n S_i ' A_{11} (S_i')^* \textrm{
and } \sum_{j=1}^p T_j A_{22} T_j ^* = \sum_{j=1}^p T_j' A_{22}
(T_j')^*.$$ Finally, for any $A_{12} \in B(H_2,H_1)$, we use our above computations
to find that
\begin{eqnarray*} \sum_{i,j} c_{ij} S_i A_{12} T_j^* & = & \sum_{i,j,k,
\ell} c_{ij}u_{ki}\overline{w_{\ell j}} S_k' A_{12} (T_\ell')^* =
\sum_{k, \ell} \Big(\sum_{i,j} c_{ij}u_{ki}\overline{w_{\ell j}}
S_k' A_{12} (T_\ell')^*\Big) \\
& = & \sum_{(k \leq p), \ell} \Big(\sum_{i,j}
c_{ij}u_{ki}\overline{w_{\ell j}} S_k' A_{12} (T_\ell')^*\Big) \\ & \ & \
+ \sum_{(k>p) , \ell} \Big(\sum_{i,j}
c_{ij}u_{ki}\overline{w_{\ell j}} S_k' A (T_\ell')^*\Big) \\
& = & \sum_{k \leq p} d_{kk} S_k' A_{12} (T_k')^* + 0 = \sum_{k=1}^p
d_{kk} S_k' A (T_k')^*.
\end{eqnarray*} 
We have shown that

\begin{displaymath}
L(A) =  \left(\begin{array}{cc} \sum_{i=1}^n S_i' A_{11} (S_i')^*
& \sum_{i=1}^p d_{ii} S_i' A_{12} (T_i') ^* \\
\sum_{i=1}^p \overline{d_{ii}} T_i' A_{21} (S_i')^* & \sum_{i=1}^p
T_i' A_{22} (T_i')^* \end{array}  \right) \end{displaymath} for all
\begin{displaymath} A =\left( \begin{array}{cc} A_{11} & A_{12} \\
A_{21} & A_{22} \end{array} \right) \in B(H_1 \oplus H_2).
\end{displaymath}
For each $i=1, \ldots, p$, define $Z_i: H_1 \oplus H_2 \rightarrow
K_1 \oplus K_2$ by
\begin{displaymath} Z_i = \left( \begin{array}{cc}d_{ii} S_i' & 0 \\ 0 &
T_i' \end{array} \right), \end{displaymath} so
\begin{displaymath} L(A) = \sum_{i=1}^p Z_i A Z_i^* + \sum_{i=1}^p \left( \begin{array}{cc}
(1-|d_{ii}|^2) S_i' A_{11} S_i'^* & 0 \\ 0 & 0
\end{array} \right) + \sum_{i=p+1}^n \left( \begin{array}{cc}
S_i' A_{11} S_i'^* & 0 \\ 0 & 0
\end{array} \right).
\end{displaymath} Since $||D|| \leq 1$, the line above shows
that $L$ is the sum of three normal completely positive maps and is 
thus normal and completely
positive.  Therefore, $\gamma$ is a corner
from $\alpha$ to $\beta$.  If, on the other hand, $n<p$, then the
same argument we just used shows that $\gamma^*$ is a corner from $\beta$ to
$\alpha$, which is equivalent to showing that $\gamma$ is a corner
from $\alpha$ to $\beta$.

For the forward direction, suppose that $\gamma$ is a corner from
$\alpha$ to $\beta$, so the map $\Upsilon: B(H_1 \oplus H_2)
\rightarrow B(K_1 \oplus K_2)$ defined by
\begin{displaymath}
\Upsilon\left( \begin{array}{cc} A_{11} & A_{12}
\\ A_{21} & A_{22} \end{array} \right)=
\left( \begin{array}{cc} \sum_{i=1}^n S_iA_{11}S_i^* &
\gamma(A_{12})
\\ \gamma^*(A_{21}) & \sum_{j=1}^p T_jA_{22}T_j^* \end{array} \right)
\end{displaymath}
is normal and completely positive.  Therefore, for some $q \in \mathbb{N}
\cup \{\infty\}$ and maps $Y_i: H_1 \oplus H_2 \rightarrow K_1
\oplus K_2$ for $i=1,2,...$, linearly independent over
$\ell_2(\mathbb{N})$, we have
$$\Upsilon(\tilde{A})= \sum_{i=1}^q Y_i\tilde{A}Y_i^*$$ for all $\tilde{A}
 \in B(H_1 \oplus H_2)$.  For $i=1,2$, let $E_i \in B(H_1 \oplus H_2)$
be projection onto $H_i$, and let $F_i \in B(K_1 \oplus K_2)$ be
projection onto $K_i$.  Since $\alpha$ and $\beta$ are contractions
we have $\Upsilon(E_1) \leq F_1$ and $\Upsilon(E_2) \leq F_2$, so
$Y_i E_j Y_i^* \leq F_j$ for each $i$ and $j$.  It follows that each
$Y_i$, $i=1, \ldots q$, can be written in the form
\begin{displaymath} Y_i  = \left( \begin{array}{cc} \tilde{S}_i & 0
\\ 0 & \tilde{T}_i \end{array} \right)\end{displaymath} for some $\tilde{S}_i \in
B(H_1, K_1)$ and $\tilde{T}_i \in B(H_2, K_2)$.

Note that $\alpha(A_{11})= \sum_{i=1}^n S_i A_{11}S_i^*= \sum_{i=1}^q \tilde{S}_i A_{11} \tilde{S_i}^*$ 
for all $A_{11}\in B(H_1)$.  For each 
$\tilde{S}_i$, define a completely positive map $L_i$ by 
$L_{i}(A) = \tilde{S}_i A\tilde{S}_i^*$
for $A \in B(H_1)$.  Since 
$\alpha - L_i$ is completely positive, it follows from the work of Arveson
in \cite{arveson} that $\tilde{S}_i$ can be written as
$$\tilde{S}_i = \sum_{j=1}^n r_{ij} S_j$$ for some complex coefficients $\{r_{ij}\}_{j=1}^n$.  
The same argument shows that for each $\tilde{T}_i$ we have
$$\tilde{T}_i = \sum_{j=1}^p b_{ij}
T_j$$ for some coefficients $\{b_{ij}\}_{j=1}^p$.  It now follows from
linear independence of the maps $\{Y_i\}_{i=1}^q$ that
$q \leq n+p$.  Let $R=(r_{ij}) \in M_{q \times n}(\C)$ and
$B=(b_{ij}) \in M_{q \times p}(\C)$, and let $A \in B(H_1)$.  We
calculate
\begin{eqnarray} \sum_{i=1}^n S_iAS_i^* &=& \sum_{i=1}^q \tilde{S}_iA\tilde{S}_i^*
= \sum_{i=1}^q \Big(\sum_{j,k=1}^n r_{ij}\overline{r_{ik}} S_j A
S_k^*
\Big) \nonumber \\
&=& \label{end?} \sum_{j,k=1}^n \Big(\sum_{i=1}^q r_{ij}
\overline{r_{ik}}\Big)S_jAS_k^*.
\end{eqnarray}
Let $M=R^T(R^T)^* \in M_{n}(\C)$, so its $jk$th entry is $ m_{jk}
= \sum_{i=1}^q r_{ij}\overline{r_{ik}}$.  Unitarily diagonalizing
$M$ as $UMU^* = D$ for some diagonal $D$ and defining maps $\{S_i'\}_{i=1}^n$ by $S_i' =
\sum_{k=1}^n \overline{u_{ik}} S_k$, we see that Eq.\eqref{end?}
and the same linear algebra technique from the proof of the backward direction
yield $$\sum_{i=1}^n S_i'AS_i'^* =
\sum_{i=1}^n S_iAS_i^* = \sum_{j,k=1}^n m_{jk} S_jAS_k^*
= \sum_{i=1}^n d_{ii} S_i'AS_i'^*.$$  Therefore $D=I$ and consequently
$M=I$, hence $||R|| =1$. An identical argument shows that $||B|| =
1$.

Let \begin{displaymath} \tilde{A} = \left( \begin{array}{cc} A_{11}
& A_{12}
\\ A_{21} & A_{22} \end{array} \right) \in B(H_1 \oplus H_2)
\end{displaymath} be arbitrary.  Let $C= (c_{jk}) \in
M_{n \times p}(\C)$ be the matrix $C=(B^*R)^T$, noting that $||C||
\leq 1$. A straightforward computation of
$\Upsilon(\tilde{A})=\sum_{i=1}^q Y_i \tilde{A} Y_i^*$ yields

\begin{eqnarray*} \gamma(A_{12})&=& \sum_{i=1}^q
S_i'A_{12}T_i'^* = \sum_{i=1}^q \Big(\Big(\sum_{j=1}^n a_{ij}
S_j\Big)A_{12}\Big(\sum_{k=1}^p \overline{b_{ik}} T_k^*\Big)\Big) \\
&=& \sum_{j,k} \Big(\sum_{i=1}^q a_{ij} \overline{b_{ik}}\Big) S_j
A_{12} T_k^* = \sum_{j,k} c_{jk} S_jA_{12}T_k^*, \end{eqnarray*}
hence $\gamma$ is of the form claimed.  \qed \end{pf}

\section{Comparison theory for $q$-positive maps}
Just as in the general study of various classes of linear operators,
it is natural to impose, and examine, an order structure for
$q$-positive maps.  If $\phi$ and $\psi$ are $q$-positive maps
acting on $M_n(\C)$, we say that $\phi$ $q$-dominates $\psi$ (and
write $\phi \geq_q \psi$) if $\phi(I + t \phi)^{-1} - \psi(I + t
\psi)^{-1}$ is completely positive for all $t \geq 0$.  We would
like to find the $q$-positive maps with the least complicated
structure of $q$-subordinates.  That last statement is not as simple
as it seems.  We might think to define a $q$-positive map $\phi$ to
be ``$q$-pure'' if $\phi \geq_q \psi \geq_q 0$ implies $\psi =
\lambda \phi$ for some $\lambda \in [0,1]$, but there exist
$q$-positive maps $\phi$ such that for every $\lambda \in (0,1)$ we
have $\phi \ngeq _q \lambda \phi$.  One such example is the Schur
map $\phi$ on $M_2 (\C)$ given by
\begin{displaymath} \phi\left( \begin{array}{cc} a_{11} &  a_{12}  \\
a_{21} & a_{22}  \\
\end{array} \right)=
\left( \begin{array}{cc} a_{11} &  (\frac{1+i}{2}) a_{12}  \\
(\frac{1-i}{2}) a_{21} & a_{22}  \\
\end{array} \right).
\end{displaymath}
As it turns out, every $q$-positive map is guaranteed to have a
one-parameter family of $q$-subordinates of a particular form:

\begin{prop}\label{qsubs}  Let $\phi \geq_q 0$.  For each $s \geq 0$,
let $\phi^{(s)}= \phi(I + s \phi)^{-1}$. Then $\phi^{(s)} \geq_q 0$
for all $s \geq 0$.  Furthermore, the set $\{\phi^{(s)}\}_{s \geq
0}$ is a monotonically decreasing family of $q$-subordinates of
$\phi$, in the sense that $\phi^{(s_1)} \geq_q \phi^{(s_2)}$ if $s_1
\leq s_2$.
\end{prop}

\begin{pf}  For all $s \geq 0$ and $t \geq 0$, we have
\begin{eqnarray*} \phi^{(s)}(I + t \phi^{(s)})^{-1} & = &
\phi(I+s \phi)^{-1}\Big(I + t \phi(I + s\phi)^{-1}\Big)^{-1} \\ &=&
\phi\Big[\Big(I + t \phi(I + s \phi)\Big)(I + s\phi)\Big]^{-1}\\ & = & \phi(I + (s+t)\phi)^{-1},
\end{eqnarray*} which is completely positive by $q$-positivity of
$\phi$.  Therefore, $\phi^{(s)} \geq_q 0$ for all $s \geq 0$.

To prove that $\phi^{(s_1)} \geq_q \phi^{(s_2)}$ if $s_1 \leq s_2$,
we let $t \geq 0$ be arbitrary and examine the map $$\Phi:
=\phi^{(s_1)}(I + t \phi^{(s_1)})^{-1} - \phi^{(s_2)}(I + t
\phi^{(s_2)})^{-1}.$$  Letting $t_1 = s_1+t$ and $t_2=s_2+t$, we
make the following observations:
\begin{equation}\label{darn1} \phi^{(s_j)}(I + t \phi^{(s_j)})^{-1}=
\phi^{(t_j)} \textrm{ for } j=1,2,
\end{equation}
\begin{equation} \label{darn2} \phi^{(t_1)} - \phi^{(t_2)} =
(I + t_2 \phi)^{-1} \Big( (I + t_2 \phi)\phi - \phi (I + t_1 \phi)
\Big) (I + t_1 \phi)^{-1}.\end{equation}

Equations \eqref{darn1} and \eqref{darn2} give us
\begin{eqnarray*} \Phi & = & (I + t_2 \phi)^{-1} \Big(
(I + t_2 \phi)\phi - \phi (I + t_1 \phi) \Big) (I + t_1 \phi)^{-1}\\
& = &
(I + t_2 \phi)^{-1} \Big((t_2-t_1) \phi^2\Big)(I + t_1 \phi)^{-1} \\
& = & (t_2-t_1) \Big(\phi(I + t_2 \phi)^{-1}\Big) \Big( \phi(I + t_1
\phi)^{-1}\Big).
\end{eqnarray*}
The last line is a non-negative multiple of a composition of
completely positive maps and is thus completely positive.  We
conclude that $\phi^{(s_1)} \geq_q \phi^{(s_2)}$.  \qed
\end{pf}

We now have the correct notion of what it means to be $q$-pure:

\begin{defn}\label{qpure}  Let $\phi: M_n(\C) \rightarrow M_n(\C)$
be unital and $q$-positive.  We say that $\phi$ is \emph{q-pure} if
its set of $q$-subordinates is precisely $\{0\} \cup
\{\phi^{(s)}\}_{s \geq 0}$.
\end{defn}

\begin{lem}\label{bdryrep2}  Let $\nu$ be a normalized unbounded boundary
weight over $L^2(0, \infty)$ of the form $$\nu(\sqrt{I - \Lambda(1)}
B \sqrt{I - \Lambda(1)}) = (f,Bf).$$  Let $\phi: M_n(\C) \rightarrow
M_n(\C)$ be a $q$-positive contraction such that $||\nu_t(I) \phi(I +
\nu_t(\Lambda(1)) \phi)^{-1}|| \leq 1$ for all
$t>0$, and let $\alpha$ be the $CP$-flow derived from the boundary
weight double $(\phi, \nu)$, with boundary generalized representation
$\pi=\{\pi_t^\#\}_{t>0}$.

Let $\beta$ be any $CP$-flow over $\C^n$, with generalized boundary
representation $\xi^\#=\{\xi_t^\#\}_{t>0}$ and boundary weight map
$\rho \rightarrow \eta(\rho)$. Then $\alpha \geq \beta$ if and only
if $\beta$ is induced by the boundary weight double $(\psi, \nu)$,
where $\psi: M_n(\C) \rightarrow M_n(\C)$ is a $q$-positive map satisfying $\phi \geq_q \psi$.
\end{lem}

\begin{pf}  As before, for each $t>0$ we let $s_t =
\nu_t(\Lambda(1))$.  Assume the hypotheses of the backward
direction.  Then $\xi_t^\#= \psi(I + s_t \psi)^{-1} \Omega_{\nu_t}$,
and the direction now follows from Theorem \ref{bdryrep} since the
line below is completely positive for all $t>0$:
$$\pi_t^\# - \xi_t^\# = (\phi(I + s_t \phi)^{-1}-\psi(I + s_t \psi)^{-1})
\Omega_{\nu_t}.$$  Now assume the hypotheses of the forward
direction.  Recall that by construction of $\nu$, the set
$\{s_t\}_{t
> 0}$ is decreasing.  If $s_t>0$ for all $t>0$ we define $P=\infty$.
Otherwise, we define $P$ to be the smallest positive number such
that $s_P = 0$.  Fix any $t_0 \in (0,P)$. Notationally, write each
$g \in H:= \C^n \otimes L^2(0, \infty)$ in its components as $g(x)=(g_1(x), \ldots, g_n(x))$, and
write $f_{t_0}$ for the function $V_{t_0}V_{t_0}^*f \in L^2(0, \infty)$, where
$V_{t_0}$ is the right shift $t_0$ units on $L^2(0, \infty)$.  Let $U_{t_0}$
be the right shift $t_0$ units on $H$.  Under our
identifications, $U_{t_0}U_{t_0}^*$ is the diagonal matrix in
$M_n(B(L^2(0, \infty)))$ with $ii^{th}$ entries $V_{t_0}V_{t_0}^*$.  Define $S: H
\rightarrow \C^n$ by
$$Sg = ((f_{t_0},g_1), \ldots, (f_{t_0}, g_n)),$$ noting
that $\Omega_{\nu_{t_0}}(A) = SAS^*$ for all $A \in B(H)$.  Since
$\phi(I + s_{t_0}\phi)^{-1}$ is completely positive, we know it has 
the form $\phi(I
+ s_{t_0} \phi)^{-1}(M) = \sum_{i=1}^m R_i M R_i^*$ for some $R_1,
\ldots, R_m \in M_n(\C)$. Therefore,
$$\pi^\#_{t_0}(A) = \Big(\phi(I + s_{t_0} \phi)^{-1}\Big) (\Omega_{\nu_{t_0}}(A))=
\sum_{i=1}^m R_iSAS^*R_i^*.$$  The map $\xi^\#_{t_0}$ is a
subordinate of $\pi^\#_{t_0}$, so from Arveson's work in metric
operator spaces in \cite{arveson}, we know that $\xi^\#_{t_0}$ has the form
$$\xi^\#_{t_0}(A)= \sum_{i,j =1}^m c_{ij} R_iSAS^*R_j^*,$$ for some
complex numbers $\{c_{ij}\}$.  Let $L_{t_0}$ be the map $L_{t_0}(M)= \sum_{i,j}
c_{ij} R_i M R_j^*$, noting that $\xi^\#_{t_0}(A)=
L_{t_0}(SAS^*)=L_{t_0}(\Omega_{\nu_{t_0}}(A))$ for all $A \in B(H)$.

Defining $\psi_{t_0}: M_n(\C) \rightarrow M_n(\C)$ by $\psi_{t_0} =
(I - \xi^\#_{t_0} \Lambda)^{-1}L_{t_0},$ we find that for arbitrary $A \in
B(H)$ and $\acute{A} \in M_n(\C)$,

\begin{eqnarray}\label{etaone} \eta_{t_0}(\rho)(A) & = & \Big(\hat{\xi}_{t_0}(I -
\hat{\Lambda} \hat{\xi}_{t_0})^{-1}\Big)(\rho)(A) \nonumber =
\rho\Big((I - \xi^\#_{t_0} \Lambda)^{-1}(\xi^\#_{t_0}(A))\Big) \\ &
= & \rho\Big((I - \xi^\#_{t_0} \Lambda)^{-1}L_{t_0}
(\Omega_{\nu_{t_0}}(A))\Big) =\rho(\psi_{t_0}(\Omega_{\nu_{t_0}}A))
\end{eqnarray} and
\begin{equation} \label{etatwo} \hat{\Lambda} \eta_{t_0} (\rho) (\acute{A})
=  \eta_{t_0}(\rho)(\Lambda(\acute{A})) =
\rho\Big(\psi_{t_0}(\Omega_{\nu_{t_0}}(\Lambda(\acute{A})))\Big) =
s_{t_0} \rho(\psi_{t_0}(\acute{A})),\end{equation} so $\hat{\Lambda}
\eta_{t_0}= s_{t_0} \hat{\psi}_{t_0}.$

Using formulas \eqref{etaone} and \eqref{etatwo} and the fact that
$\hat{\xi}_{t_0} = \eta_{t_0}(I + \hat{\Lambda} \eta_{t_0})^{-1}$,
we find
\begin{eqnarray*}\rho(\xi^\#_{t_0}) & = & \hat{\xi}_{t_0}(\rho) = \eta_{t_0}(I +
\hat{\Lambda} \eta_{t_0})^{-1}(\rho) =\Big((I + \hat{\Lambda}
\eta_{t_0})^{-1}(\rho)\Big) (\psi_{t_0} \Omega_{\nu_{t_0}})
\\ & = & \Big((I + s_{t_0} \hat{\psi}_{t_0})^{-1}(\rho)\Big)
(\psi_{t_0}\Omega_{\nu_{t_0}}) = \rho \Big((I + s_{t_0}
\psi_{t_0})^{-1} \psi_{t_0} \Omega_{\nu_{t_0}}\Big) \\ & = & \rho
\Big(\psi_{t_0}(I + s_{t_0} \psi_{t_0})^{-1} \Omega_{\nu_{t_0}}\Big)
\end{eqnarray*} for all $\rho \in M_n(\C)^*$, hence $\xi_{t_0}^\# =
\psi_{t_0}(I + s_{t_0} \psi_{t_0})^{-1}\Omega_{\nu_{t_0}}$.

We now show that the maps $\{\psi_t\}_{t>0}$ are constant on the
interval $(0, P)$. Let $t \in [t_0,P)$ be arbitrary.  For each
$\acute{A}=(a_{ij}) \in M_n(\C)$, let $A \in B(H)$ be the matrix 
with $ij$th entry $(a_{ij}/\nu_t(I))
V_{t} V_{t}^*$. Let $\rho \in M_n(\C)^*$. Straightforward
computations using formula \eqref{truncate} yield $\Omega_{t_0}(A)=
\Omega_{t}(A) = \acute{A}$ and
$\eta_{t_0}(\rho)(A)=\eta_{t}(\rho)(A)$.  Combining these equalities
gives us \begin{eqnarray*} \rho(\psi_{t_0}(\acute{A}))& = &
\rho(\psi_{t_0}\Omega_{\nu_{t_0}}(A))= \eta_{t_0}(\rho)(A)
\\ & = & \eta_{t}(\rho)(A)= \rho(\psi_{t}\Omega_{\nu_{t}}(A)) =
\rho(\psi_{t}(\acute{A})).\end{eqnarray*} Since the above formula
holds for every $\acute{A} \in M_n(\C)$ and $\rho \in M_n(\C)^*$, we
have $\psi_{t_0}=\psi_{t}$.  But both $t_0 \in (0,P)$ and $t \in
[t_0, P)$ were chosen arbitrarily, so the previous sentence shows
that $\psi_t = \psi_{t_0}$ for all $t \in (0, P)$.

Letting $\psi = \psi_{t_0}$, we have
\begin{eqnarray}\label{xi} \xi_t^\# = \psi(I + s_t \psi)^{-1}
\Omega_{\nu_t}\end{eqnarray} for all $t \in (0, P)$.  Defining
$\kappa_t$ as in the proof of Proposition \ref{bdryweight}, we
observe that $\psi(I + s_t \psi)^{-1}= \xi_t^\# \kappa_t$ for all $t
\in (0,P)$, where the right hand side is completely positive by
hypothesis. Since every $t \in (0, \infty)$ can be written as $t =
s_{t'}$ for some $t' \in (0,P)$, it follows that $\psi(I + t
\psi)^{-1}$ is completely positive for all $t>0$. Furthermore,
$\psi(I + s_t \psi)^{-1} \rightarrow \psi$ in norm as $t\rightarrow
\infty$, hence $\psi \geq_q 0$. Similarly, since $\pi_t^\#-
\xi_t^\#$ is completely positive for all $t>0$ by assumption, it
follows from our formula
$$\phi(I + s_t \phi)^{-1} - \psi(I + s_t \psi)^{-1}= (\pi_t^\# - \xi_t^
\#) \kappa_t$$ that $\phi(I + s_t \phi)^{-1} - \psi(I + s_t
\psi)^{-1}$ is completely positive for all $t>0$, and so its norm
limit (as $t \rightarrow \infty$) $\phi - \psi$ is completely
positive. Therefore, $\phi \geq_q \psi$.  Finally, since the
$CP$-flow $\beta$ is entirely determined by its generalized boundary
representation $\xi^\#$, which itself is determined by any
sequence $\{\xi_{t_n}^\#\}$ with $t_n$ tending to $0$ (see the remarks
preceding Theorem 4.29 of \cite{hugepaper}), it follows
from \eqref{xi} that $\beta$ is induced by the boundary weight
double $(\psi, \nu)$. \qed
\end{pf}

In a manner analogous to that used by Powers in \cite{bigpaper} and
\cite{hugepaper}, we define the terms \textit{$q$-corner} and
\textit{hyper maximal $q$-corner}:

\begin{defn}\label{hyp} Let $\phi: M_n(\C) \rightarrow M_n(\C)$ and $\psi:
M_k(\C) \rightarrow M_k(\C)$ be $q$-positive maps. A corner $\gamma:
M_{n \times k}(\C) \rightarrow M_{n \times k}(\C)$ from $\phi$ to
$\psi$ is said to be a $q$-corner from $\phi$ to $\psi$ if the map
\begin{displaymath} \Upsilon \left( \begin{array}{cc} A_{n \times n} &
B_{n \times k} \\ C_{k \times n} & D_{k \times k}
\end{array} \right)  = \left( \begin{array}{cc} \phi(A_{n \times n}) &
\gamma (B_{n \times k})  \\
\gamma^* (C_{k \times n}) & \psi(D_{k \times k}) \\
\end{array} \right)
\end{displaymath}
is $q$-positive.  A $q$-corner $\gamma$ is called hyper maximal if,
whenever
\begin{displaymath} \Upsilon \geq_q \Upsilon'  = \left( \begin{array}{cc} \phi' &  \gamma  \\
\gamma^* & \psi'\\
\end{array} \right) \geq_q 0,
\end{displaymath}
we have $\Upsilon = \Upsilon'$.
\end{defn}

\begin{prop}\label{basischange}  For any $q$-positive $\phi: M_n(\C)
\rightarrow M_n(\C)$ and unitary $U \in M_n(\C)$, define a map
$\phi_U$ by
$$\phi_U (A)= U^* \phi(UAU^*)U.$$
\begin{enumerate} \item The map $\phi_U$ is $q$-positive, and there
is an order isomorphism between $q$-positive maps $\beta$ such that
$\phi \geq_q \beta$ and $q$-positive maps $\beta_U$ such that
$\phi_U \geq \beta_U$.  In particular, $\phi$ is $q$-pure if and
only if $\phi_U$ is $q$-pure.
\item If $\phi$ is unital and $q$-pure, then there is a hyper maximal
$q$-corner from $\phi$ to $\phi_U$.
\end{enumerate}
\end{prop}
\begin{pf}  To prove the first assertion, we define a completely
positive map $\zeta$ on $M_n(\C)$ by $\zeta(A)= U^*AU$, noting that
$\zeta^{-1}$ is also completely positive.  For every $t \geq 0$ and
$A \in M_n(\C)$, we find that $(I + t \phi_U)^{-1}(A) = U^* (I + t
\phi)^{-1}(UAU^*)U$ and
\begin{eqnarray}\label{qposo} \phi_U(I + t \phi_U)^{-1}(A) & = & U^* \phi\Big(U(U^*(I + t
\phi)^{-1}(UAU^*)U)U^* \Big)U \nonumber \\ & = & U^*\phi(I + t
\phi)^{-1}(UAU^*)U \nonumber \\ & = & \zeta \circ \phi(I + t
\phi)^{-1} \circ \zeta^{-1}(A),
\end{eqnarray} so $\phi_U \geq_q 0$.  Given any $q$-positive map
$\beta$ such that $\phi \geq_q \beta$, define $\beta_U$ by
$\beta_U(A)= U^* \beta(UAU^*)U$.  Then $\beta_U$ is $q$-positive by
\eqref{qposo}, and for each $t \geq 0$ we have

$$\phi_U (I + t \phi_U)^{-1} - \beta_U(I + t \beta_U)^{-1} =  \zeta
\circ (\phi(I + t \phi)^{-1} - \beta(I + t \beta)^{-1}) \circ
\zeta^{-1},$$ hence $\phi_U \geq_q \beta_U$.  Of course, since
$\phi=(\phi_U)_{U^*}$, the argument just used gives an identical
correspondence between $q$-subordinates $\alpha$ of $\phi_U$ and
$q$-subordinates $\alpha_{U^*}$ of $\phi$. Our first assertion now
follows.

To prove the second statement, we define $\gamma: M_n(\C)
\rightarrow M_n(\C)$ by $\gamma(A)= \phi(AU^*)U$.  By Lemma
\ref{corners}, $\gamma$ is a corner from $\phi$ to $\phi_U$, so the
map
\begin{displaymath} \Theta \left(
\begin{array}{cc} A_{11} & A_{12} \\ A_{21} & A_{22}
\end{array} \right) = \left( \begin{array}{cc} \phi(A_{11}) &
\gamma(A_{12}) \\ \gamma^*(A_{21}) & \phi_U(A_{22})
\end{array} \right) \end{displaymath} is completely positive.  We
calculate $\gamma(I + t \gamma)^{-1}(A) = \phi(I + t
\phi)^{-1}(AU^*)U$, so for each $t \geq 0$ and $\tilde{A}=(A_{ij})
\in M_{2n}(\C)$, we have

\begin{displaymath} \Theta(I + t \Theta)^{-1} (\tilde{A})
= \left( \begin{array}{cc} \phi(I + t \phi)^{-1}(A_{11}) & \phi(I +
t \phi)^{-1}(A_{12}U^*)U \\ U^*\phi(I + t \phi)^{-1}(UA_{21}) &
\phi_U(I + t \phi_U)^{-1}(A_{22})
\end{array} \right). \end{displaymath}  This shows that $\gamma(I + t
\gamma)^{-1}$ is a corner from $\phi(I + t\phi)^{-1}$ to $\phi_U (I
+ t\phi_U)^{-1}$ for all $t \geq 0$, so $\gamma$ is a $q$-corner.
Finally, if
\begin{displaymath} \Theta' \left(
\begin{array}{cc} A_{11} & A_{12} \\ A_{21} & A_{22}
\end{array} \right) = \left( \begin{array}{cc} \alpha(A_{11}) &
\gamma(A_{12}) \\ \gamma^*(A_{21}) & \beta(A_{22})
\end{array} \right) \end{displaymath} is $q$-positive and $\Theta \geq_q 
\Theta'$, then since $\phi$ and $\phi_U$ are $q$-pure
we have $\alpha= \phi(I + t \phi)^{-1}$ for some $t \geq 0$ and
$\beta= \phi_U(I + s \phi_U)^{-1}$ for some $s \geq 0$.  Complete
positivity of $\Theta'$ implies that

\begin{displaymath} \Theta' \left(
\begin{array}{cc} I & U \\ U^* & I
\end{array} \right) = \left( \begin{array}{cc} \frac{1}{1+t}I &
U \\ U^* & \frac{1}{1+s}I
\end{array} \right) \geq 0, \end{displaymath}  so $s=t=0$ and
$\Theta=\Theta'$, hence $\gamma$ is hyper maximal. \qed \end{pf}

We have arrived at the key result of the section, which
tells us that, under certain conditions, the problem of determining
whether two $E_0$-semigroups induced by boundary weight doubles are
cocycle conjugate can be reduced to the much simpler problem of
finding hyper maximal $q$-corners between $q$-positive maps:

\begin{prop}\label{hypqc}
Let $\nu$ be a normalized unbounded boundary weight over $L^2(0,
\infty)$ which has the form $\nu(\sqrt{I - \Lambda(1)}B\sqrt{I -
\Lambda(1)})=(f, Bf)$. Let $\phi$ and $\psi$ be unital $q$-positive
maps on $M_n(\C)$ and $M_k(\C)$, respectively, and induce $CP$-flows
$\alpha$ and $\beta$ through the boundary weight doubles $(\phi,
\nu)$ and $(\psi, \nu)$.

Then $\alpha^d$ and $\beta^d$ are cocycle conjugate if and only if
there is a hyper maximal $q$-corner from $\phi$ to $\psi$.
\end{prop}

\begin{pf} Let $N=n+k$.  For the forward direction, suppose
$\alpha^d$ and $\beta^d$ are cocycle conjugate.  Since $\alpha^d$
and $\beta^d$ are of type II$_0$, we know from Theorem 
\ref{hyperflowcorn} that there is a
hyper maximal flow corner $\sigma$ from $\alpha$ to $\beta$, with associated
$CP$-flow
\begin{displaymath}\Theta =
\left (\begin{array}{cc} \alpha & \sigma
\\ \sigma^* & \beta \end{array} \right).
\end{displaymath}
Let $\Pi^\#=\{\Pi_t^\#\}$, $\pi^\#=\{\pi_t^\#\}$, and $\xi^\#=\{\xi_t^\#\}$
 be the generalized boundary
representations for $\Theta$, $\alpha$, and $\beta$, respectively. Define $s_t= \nu_t(\Lambda(1))$ for
all positive $t$, so for each $t
>0$ there is some $\mathfrak{Z}_t$ such that

\begin{displaymath}\Pi_t^\# =
\left (\begin{array}{cc} \pi_t^\# & \mathfrak{Z}_t
\\ \mathfrak{Z}_t^* & \xi_t^\# \end{array} \right) = \left(
\begin{array}{cc} \phi(I + s_t \phi)^{-1} \circ \Omega_{\nu_t, n \times n}
 & \mathfrak{Z}_t \\ \mathfrak{Z}_t^* & \psi(I + s_t \psi)^{-1}
\circ \Omega_{\nu_t, k \times k} \end{array} \right).
\end{displaymath}
\\
Since each $\mathfrak{Z}_t$ is a corner from $\phi(I + s_t
\phi)^{-1} \circ \Omega_{\nu_t, n \times n}$ to $\psi(I + s_t
\phi)^{-1} \circ \Omega_{\nu_t, k \times k}$, we have
$\mathfrak{Z}_t = L_t \circ \Omega_{\nu_t, n \times k}$ for some
$L_t$.  Define $B_t$ for each $t>0$ by
\begin{displaymath}B_t =
\left(\begin{array}{cc} \phi(I + s_t \phi)^{-1} & L_t
\\ L_t ^* & \psi(I + s_t \psi)^{-1} \end{array} \right).
\end{displaymath} We observe that $\Pi_t^\# = B_t \circ \Omega_{\nu_t,N \times N}$ for
all $t>0$, whereby the same argument given in the proof of
Lemma \ref{bdryrep2} shows that each $B_t$ has the form $B_t =
W_t (I + s_t W_t)^{-1}$ for some $W_t: M_n(\C)
\rightarrow M_n(\C)$
and that the maps $W_t$ are independent of $t$. 
Therefore, for some $\gamma: M_{n \times k}(\C) \rightarrow M_{n \times k}(\C)$,
we have
$$\mathfrak{Z}_t = \gamma(I + s_t \gamma)^{-1} \circ
\Omega_{\nu_t, n \times k}$$ for all $t>0$.  Define
$\kappa_{t, N \times N}: M_N(\C) \rightarrow B(H)$ as in Proposition
\ref{bdryweight}.  Letting

\begin{displaymath}\vartheta =
\left(\begin{array}{cc} \phi & \gamma
\\ \gamma^* & \psi \end{array} \right),
\end{displaymath}
we observe for each $t$  that $\vartheta(I+s_t \vartheta)^{-1} =
\Pi_{t}^\# \circ \kappa_{t, N \times N}$ is the composition of
completely positive maps and is thus completely positive, hence
$\vartheta \geq_q 0$.  Suppose that for some map $\vartheta'$ we
have

\begin{displaymath}\vartheta \geq_q \vartheta' =
\left (\begin{array}{cc} \phi' & \gamma
\\ \gamma^* & \psi' \end{array} \right) \geq_q 0.
\end{displaymath}
As in Proposition \ref{bdryweight}, the boundary weight map $\rho
\in M_{N}(\C)^* \rightarrow L(\rho)$ defined by $L(\rho)(C)=
\rho(\vartheta'(\Omega_{\nu,N \times N} (C))$ induces a $CP$-flow
$\Theta'$ over $\C^N$, where for some $CP$-flows $\alpha'$ over $\C^n$ and $\beta'$
over $\C^k$, we have

\begin{displaymath}\Theta' =
\left (\begin{array}{cc} \alpha' & \sigma
\\ \sigma^* & \beta' \end{array} \right).
\end{displaymath}
By Lemma \ref{bdryrep2}, we have $\Theta \geq \Theta'$
since $\vartheta \geq_q \vartheta'$.
But $\Theta$ is a hyper maximal flow corner, so
$\Theta=\Theta'$.  Our formulas for the generalized boundary representations
imply that $\phi(I + t \phi)^{-1}= \phi'(I + t \phi')^{-1}$ 
and $\psi(I + t \psi)^{-1} = \psi'(I + t \psi')^{-1}$ for all $t>0$,
hence $\phi=\phi'$ and $\psi=\psi'$.  We
conclude that $\gamma$ is a hyper maximal $q$-corner.\\ \\
For the backward direction, suppose there is a hyper maximal
$q$-corner $\gamma$ from $\phi$ to $\psi$, so the map $\Upsilon:
M_{N}(\C) \rightarrow M_{N}(\C)$ defined by
\begin{displaymath} \Upsilon \left( \begin{array}{cc} A_{n \times n} &
B_{n \times k} \\ C_{k \times n} & D_{k \times k}
\end{array} \right)  = \left( \begin{array}{cc} \phi(A_{n \times n}) &
\gamma (B_{n \times k})  \\
\gamma^* (C_{k \times n}) & \psi(D_{k \times k}) \\
\end{array} \right)
\end{displaymath} is $q$-positive.  By Proposition \ref{bdryweight},
the boundary weight map $\rho \in M_N(\C)^* \rightarrow \Xi(\rho)$ 
defined by
$$\Xi(\rho)(A) = \rho(\Upsilon(\Omega_{\nu, N \times N}(A)))$$ is the boundary weight map of a $CP$-flow
$\theta$ over $\C^{N}$, where for some $\Sigma$ we have

\begin{displaymath}\theta =
\left (\begin{array}{cc} \alpha & \Sigma
\\ \Sigma^* & \beta \end{array} \right).
\end{displaymath}

Let \begin{displaymath}\theta' = \left (\begin{array}{cc} \alpha' &
\Sigma
\\ \Sigma^* & \beta' \end{array} \right)
\end{displaymath}
be any $CP$-flow such that $\theta \geq \theta'$. Letting
$\mathcal{Z}_t= \gamma(I + s_t \gamma)^{-1} \circ \Omega_{\nu_t,n
\times k}$ for all $t>0$, we see the generalized boundary representations
$\Pi^\#=\{\Pi_t ^\#\}$ and $\Pi'=\{\Pi_t'\}$ for $\theta$ and $\theta'$
satisfy

\begin{displaymath}\Pi_t ^\# = \left (\begin{array}{cc} \pi_t^\# &
\mathcal{Z}_t
\\ \mathcal{Z}_t^* & \xi_t ^\# \end{array} \right)
\geq \Pi_t' = \left (\begin{array}{cc} \pi_t' & \mathcal{Z}_t
\\ \mathcal{Z}_t^* & \xi_t' \end{array} \right)
\end{displaymath}
for all $t> 0$.  Lemma \ref{bdryrep2} implies that for some $\phi'$
and $\psi'$ with $\phi \geq_q \phi' \geq_q 0$ and $\psi \geq_q \psi'
\geq_q 0$ we have $\pi_t' = \phi'(I + s_t \phi')^{-1}\circ
\Omega_{\nu_t, n \times n}$ and $\xi_t' = \psi'(I + s_t
\psi')^{-1}\circ \Omega_{\nu_t, k \times k}$ for all $t>0$.  Defining
$\Upsilon': M_N(\C) \rightarrow M_N(\C)$ by 
\begin{displaymath} \Upsilon' \left( \begin{array}{cc} A_{n \times n} &
B_{n \times k} \\ C_{k \times n} & D_{k \times k}
\end{array} \right)  = \left( \begin{array}{cc} \phi'(A_{n \times n}) &
\gamma (B_{n \times k})  \\
\gamma^* (C_{k \times n}) & \psi'(D_{k \times k}) \\
\end{array} \right),
\end{displaymath}
we observe that $\Pi_{t}' \circ \kappa_{\nu_t,N \times N}= \Upsilon'(I +
s_t \Upsilon')^{-1}$ for all $s_t>0$, hence $\gamma$ is a $q$-corner from
$\phi'$ to $\psi'$.  Hyper maximality of $\gamma$ implies $\phi=\phi'$
and $\psi=\psi'$, thus $\theta=\theta'$. Therefore, $\sigma$ is a
hyper maximal flow corner from $\alpha$ to $\beta$, so $\alpha^d$ and
$\beta^d$ are cocycle conjugate by Theorem \ref{hyperflowcorn}.
\qed \end{pf}

\section{$E_0$-semigroups obtained from rank one unital $q$-pure maps}
Any unital linear map $\phi: M_n(\C) \rightarrow M_n(\C)$ of rank
one is of the form $\phi(A) = \tau(A)I$ for some linear functional
$\tau$.  If $\phi$ is positive, then $\tau$ is positive and
$\tau(I)=1$, so $\tau$ is a state.  On the other hand, given any
state $\rho$, the map $\phi$ defined by $\phi(A)=\rho(A)I$ is unital
and completely positive. Furthermore, $\phi$ is $q$-positive since
$\phi(I+t \phi)^{-1}= (1/(1+t)) \phi$ for all $t>0$.  The rank one unital
$q$-positive maps are therefore precisely the maps $A
\rightarrow \rho(A) I$ for states $\rho$.

The goal of this section is to determine when such maps are
$q$-pure, and then to determine when the
$E_0$-semigroups induced by $(\phi, \nu)$ and $(\psi, \nu)$ are cocycle conjugate, 
where $\phi$ and $\psi$ are rank one unital $q$-pure maps and
$\nu$ is a normalized unbounded boundary weight of the form
$\nu(\sqrt{I - \Lambda(1)} B \sqrt{I - \Lambda(1)}) = (f,Bf)$ (Theorem
\ref{statesbig}). We also 
obtain a partial result for comparing $E_0$-semigroups induced by $(\phi, \nu)$
and $(\psi, \mu)$ for rank one unital $q$-pure maps $\phi$ and $\psi$ and any 
normalized unbounded boundary
weights $\nu$ and $\mu$ over $L^2(0, \infty)$ (Corollary \ref{thesearenew}).  \\
\\ We begin with a lemma:
\begin{lem}\label{tauca} Let $\rho$ be a faithful state on $M_n(\C)$, and
define a unital $q$-positive map $\phi:M_n(\C) \rightarrow M_n(\C)$
by $\phi(A)= \rho(A) I$. For any non-zero positive linear functional
$\tau$ on $M_n(\C)$ and non-zero positive operator $C \in M_n(\C)$,
define $\psi_{\tau, C}: M_n(\C) \rightarrow M_n(\C)$ by $\psi_{\tau,
C}(A) = \tau(A)C$.

Then $\psi_{\tau, C}$ is $q$-positive, and $\phi \geq_q \psi_{\tau,
C}$ if and only if $\psi_{\tau, C} = \lambda \phi$ for some $\lambda
\in (0,1]$.
\end{lem}

\begin{pf} Note that for all $A \in M_n(\C)$ and $t
\geq 0$, we have $(I+t \psi_{\tau, C})^{-1}(A)= A - t
\tau(A)/(1+t \tau(C))C$, so \begin{eqnarray}\label{tauc}\psi_{\tau,
C}(I+t \psi_{\tau, C})^{-1}(A) = \frac{\tau(A)}{1 + t \tau(C)}
C,\end{eqnarray} hence $\psi_{\tau, C}$ is $q$-positive.  It follows
from \eqref{tauc} that $\phi(I+t\phi)^{-1}(A) = (\rho(A)/(1+t))
I$ for all $A \in M_n(\C)$.

Assume the hypotheses of the forward direction.  Since $\phi \geq_q
\psi_{\tau, C}$, we have \begin{equation}\label{first} \frac{\rho(A)
I}{1+t} \geq \frac{\tau(A)C}{1 + t \tau (C)}\end{equation} for all
$t \geq 0$ and $A \geq 0$.  This is impossible if $\tau(C)=0$, so we
may assume $\tau(C) \neq 0$.  Letting $t \rightarrow \infty$ in
\eqref{first} yields \begin{equation} \label{second} \rho(A)I\geq
\frac{\tau(A) C}{\tau(C)}\end{equation} for all $A \geq 0$. Setting
$A=C$ in \eqref{second}, we see $\rho(C) I - C \geq 0$, yet
$$\rho \Big(\rho(C)I-C\Big)= \rho(C)-\rho(C)= 0,$$ hence
$C=\rho(C)I$ by faithfulness of $\rho$.  Rewriting \eqref{second} as $$\rho(A) I \geq
\frac{\tau(A)}{\tau(\rho(C)I)} \rho(C)I = \frac{\tau(A)}{||\tau||}I$$
for all $A\geq 0$, we see that $\rho-\tau / ||\tau||$ is a
positive linear functional. Therefore, $$\Big|\Big|\rho -
\frac{\tau}{||\tau||}\Big|\Big| = \rho(I) - \frac{\tau(I)}{||\tau||}
= 1-1=0,$$ hence $\tau = ||\tau||\rho$.   Setting $t=0$ and $A=I$ in
\eqref{first} gives us $||\tau||= \tau(I) = \lambda / \rho(C)$
for some $\lambda \in (0,1]$.  Therefore,
$$\psi_{\tau, C} (A) = \tau(A)C =||\tau|| \rho(A) \rho(C) I = \lambda
\rho(A)I = \lambda \phi(A)$$ for all $A \in M_n(\C)$, proving the
forward direction.

The backward direction follows from Proposition \ref{qsubs} since
$\lambda \phi = \phi^{(-1 + 1/\lambda)}$ for every
$\lambda \in (0,1]$.
 \qed \end{pf}

{\bf Remark:} Let $\psi: M_n(\C) \rightarrow M_n(\C)$ be a non-zero
$q$-positive
contraction such that the maps $L_{\psi_t} : = t \psi(I + t
\psi)^{-1}$ satisfy $||L_{\psi_t}|| <1$ for all $t
> 0$.  By compactness of the unit ball of $B(M_n(\C))$, the maps $L_{\psi_t}$
have some norm limit as $t \rightarrow \infty$.  This limit is
unique: Pick any orthonormal basis with respect to the trace inner
product $(A,B)= tr(A^*B)$ of $M_n(\C)$, and let $M_t$ be the $n^2
\times n^2$ matrix of $L_{\psi_t}$ with respect to this basis. From
the cofactor formula for $(I + t \psi)^{-1}$, we know that the
$ij$th entry of $M_t$ is a rational function $r_{ij}(t)$.
Uniqueness of $\lim_{t \rightarrow \infty}L_{\psi_t}$ now follows
from the fact that each $r_{ij}(t)$ has a unique limit as $t
\rightarrow \infty$. We call this limit $L_\psi$.  Noting that
$$t \psi
= L_{\psi_t}(I - L_{\psi_t})^{-1} = L_{\psi_t} + L_{\psi_t}^2 +
\ldots$$ for each $t>0$, we claim that $L_{\psi}$ fixes a positive
element $T$ of norm one.  To prove this, we first observe for each $k \in \mathbb{N}$ 
and $t>0$ that
\begin{eqnarray*}t ||\psi|| & =& t ||\psi(I)|| \leq ||L_{\psi_t}(I)|| + \ldots +
||(L_{\psi_t})^{k-1}(I)|| + k \sum_{n=1}^\infty ||(L_{\psi_t})^{kn}(I)|| \\
& < & (k-1) + k \sum_{n=1}^\infty ||(L_{\psi_t})^k(I)||^n
,\end{eqnarray*} hence
$$1= \lim_{t\rightarrow \infty} ||(L_{\psi_t})^k(I)|| =
||(L_{\psi})^k(I)||.$$  
Therefore, all elements of the sequence $\{T_k\}_{k \in
\mathbb{N}}$ defined by $T_k = (L_{\psi})^k(I)$ satisfy $T_k \geq 0$
and $||T_k||=1$.  Since $T_k-T_{k+1} = (L_{\psi})^k(I-T_1) \geq 0$
for all $k$, the sequence $\{T_k\}_{k \in \mathbb{N}}$ is
monotonically decreasing and therefore has a positive norm limit $T$
with $||T||=1$.  Finally, $L_\psi$ fixes $T$ since
$L_{\psi}(T)=\lim_{k \rightarrow \infty} L_{\psi}^{k+1}(I)=T$. The
information at hand suffices in showing that a large class of maps
is $q$-pure:

\begin{prop}\label{littlestates}  Let $\rho$ be a state on $M_n(\C)$, and define a
$q$-positive map $\phi$ on $M_n(\C)$ by $\phi(A)= \rho(A) I$.  Then
$\phi$ is $q$-pure if and only if $\rho$ is faithful.
\end{prop}

\begin{pf}  For the forward direction, we prove the
contrapositive.  If $\rho$ is not faithful, then for some $k<n$ and
mutually orthogonal vectors $f_1, \ldots, f_k$ with $\sum_{i=1}^k
||f_i||^2=1$, we have $\rho(A)= \sum_{i=1}^k (f_i, Af_i)$ for all $A
\in M_n(\C)$. Let $P$ be the projection onto the $k$-dimensional
subspace of $\C^n$ spanned by the vectors $f_1, \ldots, f_k$, and
define a $q$-positive map $\psi: M_n(\C) \rightarrow M_n(\C)$ by
$\psi(A) = \rho(A) P$.  For each $t \geq 0$ and $A \in M_n(\C)$, we
find
$$(\phi^{(t)}-\psi^{(t)})(A)= \frac{1}{1+t}(\phi(A)-\psi(A)) =
\frac{1}{1+t} \rho(A) (I-P),$$  so $\phi \geq_q \psi$. Obviously,
$\psi
\neq \phi^{(s)}$ for any $s \geq 0$, so $\phi$ is not $q$-pure. \\
\\ To prove the backward direction,
 suppose $\phi \geq_q \psi \geq_q 0$ for some $\psi \neq 0$, and form
$L_\psi$ and $L_\phi$.  Since $L_{\phi_t} = (t/(1+t))\phi$ for
each $t>0$, we have $L_\phi=\phi$.  The map $L_{\phi_t}-L_{\psi_t}$
is completely positive for all $t$, so by taking its limit as $t \rightarrow
\infty$ we see $\phi - L_\psi$ is
completely positive. By the remarks preceding this proposition, we
know that $L_\psi$ fixes a positive $T$ with $||T||=1$. But $(\phi
-L_\psi)(T) = \rho(T)I - T \geq 0$, so $\rho(T)=1$, hence $T=I$ by
faithfulness of $\rho$.

By complete positivity of $\phi-L_\psi$, we have $||\phi - L_\psi||
= ||\phi(I) - L_\psi(I)|| = 0$, so $\phi = L_{\psi}$. Therefore,
\begin{eqnarray} 0 & = & \lim_{t \rightarrow \infty} \Big((\phi -
L_{\psi_t}) \Big(\frac{I}{t} + \psi\Big)\Big) = \lim_{t \rightarrow
\infty} \Big( \phi\Big(\frac{I}{t} + \psi\Big) -
L_{\psi_t}\Big(\frac{I}{t}
+ \psi\Big) \Big) \nonumber \\
& = & \lim_{t \rightarrow \infty} \Big( \frac{\phi}{t} + \phi \psi -
t \psi(I + t \psi)^{-1}\Big(\frac{I}{t} + \psi\Big)\Big) =
\lim_{t\rightarrow \infty} \frac{\phi}{t} + \phi \psi - \psi \nonumber\\
\label{third} & = & \phi \psi - \psi.
\end{eqnarray}

Letting $\tau$ be the positive linear functional $\tau= \rho \circ
\psi$, we conclude from \eqref{third} that $\psi(A)=\rho(\psi(A)) I
= \tau(A)I$ for all $A \in M_n(\C)$. Lemma \ref{tauca} implies that
$\psi = \lambda \phi = \phi^{(-1+1/\lambda)}$ for some
$\lambda \in (0,1]$.
 \qed \end{pf}

To prove the main result of the section, we need the following:
\begin{lem}\label{generalnormone}  Let $\phi: M_n(\C) \rightarrow
M_n(\C)$ and $\psi: M_k(\C) \rightarrow M_k(\C)$ be rank one unital
$q$-pure maps, and let $\nu$ and $\mu$ be normalized
unbounded boundary weights over
$L^2(0, \infty)$. If the boundary weight doubles $(\phi, \nu)$ and
$(\psi, \mu)$ induce cocycle conjugate $E_0$-semigroups $\alpha^d$
and $\beta^d$, then there is a corner $\gamma$ from $\phi$ to $\psi$
such that $||\gamma||=1$.
\end{lem}
\begin{pf}  By construction, $\alpha^d$ and $\beta^d$ 
are type II$_0$ $E_0$-semigroups.  If they are cocycle conjugate, then
by Theorem \ref{hyperflowcorn},
there is a hyper maximal flow corner $\sigma$ from $\alpha$ to
$\beta$ with associated $CP$-flow $\Theta$ over $K_1 \oplus K_2$, where
\begin{displaymath}\Theta =
\left (\begin{array}{cc} \alpha & \sigma
\\ \sigma^* & \beta \end{array} \right).
\end{displaymath}
Let $H_1 = \C^n \otimes L^2(0, \infty)$ and $H_2 = \C^k \otimes L^2(0, \infty)$.
Write the boundary representation $\Pi=\{\Pi_t^\#\}$ for $\Theta$ as

\begin{displaymath} \Pi_t^\# = \left (\begin{array}{cc} \frac{1}{1+
\nu_t(\Lambda(1))} \phi \circ \Omega_{\nu_t, n \times n} & \mathfrak{Z}_t \\
\mathfrak{Z}_t^* & \frac{1}{1+\mu_t(\Lambda(1))}\psi \circ
\Omega_{\mu_t, k \times k}\end{array} \right)
\end{displaymath}
for some maps $\{\mathfrak{Z}_t\}_{t>0}$ from $B(H_2, H_1)$ into $B(K_2, K_1)$.  Let $\rho_{11}
\rightarrow \omega(\rho_{11})$ and $\rho_{22} \rightarrow
\eta(\rho_{22})$ denote the boundary weight maps for $\alpha$ and
$\beta$, respectively.  Let $\rho \rightarrow \Xi(\rho)$ be the
boundary weight map for $\Theta$, so for some map $\rho_{12}
\rightarrow \ell(\rho_{12})$ from $M_{n \times k}(\C)^*$ to weights
on $B(H_2, H_1)$ we have
\begin{displaymath}
\Xi \left (\begin{array}{cc} \rho_{11} & \rho_{12} \\ \rho_{21} &
\rho_{22}
\end{array} \right)= \left (\begin{array}{cc} \omega(\rho_{11}) &
\ell(\rho_{12})
\\ \ell^*(\rho_{21}) & \eta(\rho_{22})
\end{array} \right).
\end{displaymath}

Denote by $U_t$ the right shift $t$ units on $H$, and let $\pi^\#$
and $\xi^\#$ be the generalized boundary representations for $\alpha$
and $\beta$, respectively.  For every
$A=(A_{ij}) \in \mathcal \bigcup_{t>0} U_t B(H) U_t^*$ and bounded
family of functionals $\{\rho(t)=(\rho_{ij}(t))\}_{t > 0}$ in
$M_{n+k}(\C)^*$, we observe that the argument used in
Corollary 3.3 to show that $\pi_0 ^\# = \xi_0^\# = 0$ implies
$$\lim_{t \rightarrow 0}\omega_t(I +
\hat{\Lambda}\omega_t)^{-1}(\rho_{11}(t))(A_{11}) = \lim_{t
\rightarrow 0} \eta_t(I +
\hat{\Lambda}\eta_t)^{-1}(\rho_{22}(t))(A_{22})=0,$$ so by complete
positivity of the generalized boundary representation, we have
\begin{eqnarray}\label{ell} \lim_{t \rightarrow 0} \ell_t(I +
\hat{\Lambda}\ell_t)^{-1}(\rho_{12}(t))(A_{12})=0.\end{eqnarray}

We claim that $\rho_{12} \rightarrow \ell(\rho_{12})$
is unbounded.  If $\ell$ is bounded, then for each $\rho_{12} \in M_{n \times
k}(\C)^*$, the family $\rho_{12}(t):=(I +
\hat{\Lambda}\ell_t)(\rho_{12})$ is bounded, and it follows from
\eqref{ell} that
\begin{eqnarray}\label{ell22} \lim_{t \rightarrow 0}
\ell_t(\rho_{12})(A_{12})=0 \end{eqnarray} for each $A_{12} \in
\bigcup_{t>0} W_t B(H_2, H_1) X_t^*$, where $W_t$ and $X_t$ are the
right shift $t$ units on $H_1$ and $H_2$, respectively.  Let $A_{12}
\in \bigcup_{t>0} W_t B(H_2, H_1) X_t^*$, so $A_{12} = W_s B X_s^*$
for some $s>0$ and $B \in B(H_2,H_1)$.  For all $b<s$, we have
\begin{eqnarray*} \ell_b(\rho_{12})(A_{12}) & = & \ell_b(\rho_{12})(W_sBX_s^*)
= \ell(\rho_{12})(W_bW_b^*W_sBX_s^*X_bX_b^*)\\ & = &
\ell(\rho_{12})(W_b W_{s-b} B X_{s-b}^*X_b^*) =  \ell(\rho_{12})(W_s
B X_s^*) \\ &  = & \ell(\rho_{12})(A_{12}). \end{eqnarray*}
Therefore, by equation \eqref{ell22} we have
$\ell(\rho_{12})(A_{12})=0.$  Let $A \in B(H_2, H_1)$, $\rho_{12}
\in M_{n \times k}(\C)^*$, and $t>0$ be arbitrary. From above we
have $$\ell_t(\rho_{12})(A)= \ell(\rho_{12})(W_tAX_t^*)=0,$$ hence
$\ell_t \equiv 0$ for all $t>0$.  We conclude from uniqueness of the
generalized boundary representation that $\rho_{12} \rightarrow
\ell(\rho_{12})$ is the zero map.  The boundary weight map $\rho
\rightarrow \Xi'(\rho)$ defined by

\begin{displaymath}
\Xi' \left (\begin{array}{cc} \rho_{11} & \rho_{12} \\ \rho_{21} &
\rho_{22}
\end{array} \right)= \left (\begin{array}{cc} \omega(\rho_{11}) &
0
\\ 0 & 0
\end{array} \right)
\end{displaymath}
gives rise to the $CP$-flow
\begin{displaymath}\Theta' =
\left (\begin{array}{cc} \alpha & \sigma
\\ \sigma^* & \beta', \end{array} \right)
\end{displaymath}
where $\beta'$ is the non-unital $CP$-flow
$\beta_t'(A_{22})=X_tA_{22}X_t^*$. Trivially, $\Theta \neq \Theta'$
and $\Theta \geq \Theta'$, contradicting hyper maximality of
$\sigma$. Therefore, the map $\rho_{12} \rightarrow \ell(\rho_{12})$
is unbounded.

Since $\Pi_t^\#$ is a contraction for every $t>0$, so is
$\mathfrak{Z}_t$, hence the map $\mathfrak{Z}_t \circ \Lambda: M_{n
\times k}(\C) \rightarrow M_{n \times k} (\C)$ is a contraction for
each $t>0$.  A compactness argument shows that $\mathfrak{Z}_{t_n}
\circ \Lambda$ has a norm limit $\gamma$ for some sequence $\{t_n\}$
tending to zero, where $||\gamma|| \leq 1$. From unboundedness of
$\ell$ and the formula $\ell_t=\hat{\mathfrak{Z}}_t(I
-\hat{\Lambda}\hat{\mathfrak{Z}}_t)^{-1}$ for all $t>0$, it follows
that $I-\gamma$ is not invertible, so $||\gamma|| \geq 1$, hence
$||\gamma||=1$.  We claim that $\gamma$ is a corner from $\phi$ to
$\psi$.  Indeed, for the family of completely positive maps
$\{R_t\}_{t>0}$ defined by $R_t = \Pi_t^\# \circ \Lambda$, we have
\begin{displaymath} \lim_{n \rightarrow \infty} R_{t_n} = \lim_{n
\rightarrow \infty} \left( \begin{array}{cc}
\frac{\nu_{t_n}(\Lambda(1))}{1+\nu_{t_n}(\Lambda(1))} \phi & \mathfrak{Z}_{t_n} \circ \Lambda \\
(\mathfrak{Z}_{t_n} \circ \Lambda)^* &
\frac{\mu_{t_n}(\Lambda(1))}{1+\mu_{t_n}(\Lambda(1))} \psi.
\end{array} \right) = \left( \begin{array}{cc}\phi & \gamma \\
\gamma^* & \psi  \end{array} \right).
\end{displaymath}
\qed \end{pf}

If \ $\nu$ \ is \ a normalized \ unbounded \ boundary \ weight \ over \ $L^2(0, \infty)$ \
of \ the \ form $\nu(\sqrt{I - \Lambda(1)}B\sqrt{I - \Lambda(1)})=(f,Bf)$ and if $\phi: M_n(\C) \rightarrow M_n(\C)$
is unital and $q$-pure, we know from Propositions \ref{basischange}
and \ref{hypqc} 
that the condition $\psi=\phi_U$ is sufficient for the boundary weight
doubles $(\phi, \nu)$ and 
$(\psi, \nu)$ to induce cocycle conjugate $E_0$-semigroups.  In the case
that $\phi$ is a rank one unital $q$-pure map, this condition is also necessary: 
\begin{thm}\label{statesbig}  Let $\phi_1: M_n(\C) \rightarrow M_n(\C)$ and 
$\phi_2: M_k(\C) \rightarrow M_k(\C)$ be rank
one unital $q$-pure maps.  Let
$\nu$ be a normalized unbounded boundary weight over $L^2(0,
\infty)$ of the form $\nu(\sqrt{I - \Lambda(1)} B \sqrt{I -
\Lambda(1)}) = (f,Bf)$.

Then the boundary weight doubles $(\phi_1, \nu)$ and $(\phi_2, \nu)$
induce cocycle conjugate $E_0$-semigroups if and only if $n=k$ and
$\phi_2=(\phi_1)_U$ for some unitary $U \in M_n(\C)$.
\end{thm}
\begin{pf} The backward direction follows immediately from
Propositions \ref{basischange} and \ref{hypqc}.  Assume the
hypotheses of the forward direction. Since $\phi_1$ and $\phi_2$ are
rank one, unital, and $q$-pure, there exist faithful states $\rho_1$
on $M_n(\C)$ and $\rho_2$ on $M_k(\C)$ such that $\phi_1(M)=
\rho_1(M)I_{n \times n}$ and $\phi_2(B) = \rho_2(B)I_{k \times k}$
for all $M \in M_n(\C)$, $B \in M_k(\C)$. By Lemma
\ref{generalnormone}, there is a corner $\gamma$ from $\phi_1$ to
$\phi_2$ such that $||\gamma||=1$. Therefore, for some $A_0 \in M_{n
\times k}(\C)$ of norm one and unit vectors $f_0 \in \C^n$ and $g_0
\in \C^k$, we have $|(f_0, \gamma(A_0)g_0)|=1$. Define $\omega \in
M_{n \times k}(\C)^*$ by $\omega(A) = (f_0, \gamma(A)g_0)$, noting
that $||\omega||=|\omega(A_0)|=1$. We claim that the map
$\tilde{\psi}: M_{n+k}(\C) \rightarrow M_2(\C)$ defined by
\begin{displaymath} \tilde{\psi} \left( \begin{array}{cc} A_{11} & A_{12}
\\ A_{21} & A_{22} \end{array}  \right)=
\left( \begin{array}{cc} \rho_1(A_{11}) & \omega(A_{12}) \\
\omega^*(A_{21}) & \rho_2(A_{22})  \end{array} \right)
\end{displaymath} is completely positive.  To see this, let
$\{\tilde{F_i}\}_{i=1}^\ell$ be arbitrary vectors in $\C^2$, writing
each $\tilde{F_i}$ as
\begin{displaymath}\tilde{F_i}
= \left( \begin{array}{cc} \lambda_{1i}\\
\lambda_{2i} \end{array} \right) \end{displaymath} for some complex
numbers $\{\lambda_{1i}\}_{i=1}^\ell$ and $\{\lambda_{2i}
\}_{i=1}^\ell$.

Since the map $\psi: M_{n+k}(\C) \rightarrow M_{n+k}(\C)$ defined by
\begin{displaymath} \psi \left( \begin{array}{cc} A_{11} & A_{12}
\\ A_{21} & A_{22} \end{array}  \right)=
\left( \begin{array}{cc} \rho_1(A_{11})I & \gamma(A_{12}) \\
\gamma^*(A_{21}) & \rho_2(A_{22})I  \end{array} \right)
\end{displaymath} is completely positive by assumption, we know
that for any  $A_1, \ldots, A_\ell \in M_{n+k}(\C)$ and the vectors
\begin{displaymath}  F_i = \left( \begin{array}{cc} \lambda_{1i} f_0 \\
\lambda_{2i} g_0 \end{array} \right) \in \C^{n+k}, \ \ i=1, \ldots,
k,
\end{displaymath} we have
$$ \sum_{i,j=1}^\ell \Big(F_i, \psi(A_i^*A_j) F_j\Big) \geq 0.$$
However, for each $i$ and $j$ we find that \begin{eqnarray*}
\Big(F_i, \psi(A_i^*A_j)F_j\Big)_{\C^{n+k}}& = &
\overline{\lambda_{1i}}\lambda_{1j} \rho_1((A_i^*A_j)_{11}) +
\overline{\lambda_{1i}}\lambda_{2j} \omega((A_i^*A_j)_{12})\\ & \ &
+ \overline{\lambda_{2i}}\lambda_{j1}
\overline{\omega([(A_i^*A_j)_{21}] ^*)} +
\overline{\lambda_{2i}}\lambda_{2j} \rho_2((A_i^*A_j)_{22}) \\
& = & \Big(\tilde{F_i}, \tilde{\psi}(A_i^*A_j) \tilde{F_j}
\Big)_{\C^2}.
\end{eqnarray*}
Therefore, for all $\ell \in \mathbb{N}$, 
$A_1, \dots, A_\ell \in
M_{n+k}(\C)$, and $\tilde{F_1}, \ldots, \tilde{F_\ell} \in \C^{2}$
, we have $\sum_{i,j=1}^\ell \Big(\tilde{F_i}, \tilde{\psi}
(A_i^*A_j) \tilde{F_j}\Big) \geq 0$, 
so $\tilde{\psi}: M_{2n}(\C) \rightarrow M_2(\C)$ is
completely positive.  Since $\rho_1$ and $\rho_2$ are positive
linear functionals (hence completely positive maps), $\omega$ is a
corner from $\rho_1$ to $\rho_2$.

By faithfulness of $\rho_1$ and $\rho_2$, there exist monotonically
increasing sequences of strictly positive numbers
$\{\lambda_i\}_{i=1}^n$ and $\{\mu_j\}_{j=1}^k$ with $\sum_{i=1}^n
\lambda_i^2 = \sum_{j=1}^k \mu_j^2 =1$,  along with orthonormal sets
of vectors $\{f_i\}_{i=1}^n$ and $\{g_j\}_{j=1}^k$, such that
$\rho_1(M)= \sum_{i=1}^n \lambda_i^2 (f_i, M f_i)$ and $\rho_2(B) =
\sum_{j=1}^k \mu_j^2 (g_j, B g_j)$ for all $M \in M_n(\C)$, $B \in
M_k(\C)$. Given $A \in M_{n \times k}(\C)$, let $\tilde{A}$ be the
matrix whose $ji$th entry is $(f_i, A g_j)$, observing that
$||\tilde{A}|| =||A||$. Let $D_\lambda$ and $D_\mu$ be the diagonal
matrices whose $ii$th entries are $\lambda_i$ and $\mu_i$,
respectively, for all $i$, and let $D_{\lambda^2}$ and $D_{\mu^2}$
be the diagonal matrices whose $ii$th entries are $\lambda_i^2$ and
$\mu_i^2$, respectively, observing that $D_{\lambda^2} =
(D_\lambda)^2$ and $D_{\mu^2} = (D_{\mu})^2$.\\ \\
By Proposition \ref{corners}, $\omega$ has the form
$$\omega(A)= \sum_{i,j} c_{ij} \lambda_i \mu_j (f_i, A g_j) =
tr(C D_\mu \tilde{A} D_\lambda)= tr\Big(C D_\mu (D_\lambda
\tilde{A}^*)^* \Big)$$ for some $C=(c_{ij}) \in M_{n\times k}(\C)$
such that $||C|| \leq 1$.

By the Cauchy-Schwartz inequality for the inner product $(B,A)=tr(AB^*)$ on $M_{n
\times k} (\C)$, we have
\begin{eqnarray}1&=& |\omega(A_0)|^2 =  |tr(C D_\mu (D_\lambda
\tilde{A_0}^*)^*) |^2 =|(C D_\mu, D_\lambda \tilde{A_0}^*)|^2 \nonumber \\
& \leq & (C D_\mu, C D_\mu) (D_\lambda
\tilde{A_0}^*,D_\lambda \tilde{A_0}^*)= tr(D_\mu C^* C D_\mu)
tr(D_\lambda
\tilde{A_0}^*\tilde{A_0} D_\lambda) \nonumber \\
\label{ii} & \leq & tr(D_{\mu^2}I_k)tr(D_{\lambda^2}I_n) \leq 1*1=1.
\end{eqnarray}
Since equality holds in all the inequalities above, we have $m C D_\mu = D_\lambda
\tilde{A_0}^*$ for some $m \in \C$.  It follows from \eqref{ii}
that $|m|=1$ since $||CD_\mu||_{tr}=||D_\lambda \tilde{A_0}^*||_{tr}=1$.  
Furthermore, since equality holds in \eqref{ii} and the trace map is faithful, we
have $C^*C=I_k$ and $\tilde{A_0}^*\tilde{A_0}=I_n$.  But $C \in M_{n
\times k}(\C)$ and $\tilde{A_0}^* \in M_{n \times k}(\C)$, so $n=k$,
hence $C$ and $\tilde{A_0}$ are unitary. 

Writing $D_\lambda = m C D_\mu \tilde{A_0}= (m C
\tilde{A_0})(\tilde{A_0}^* D_\mu \tilde{A_0}),$ we observe that $m C \tilde{A_0}$
is unitary and $\tilde{A_0}^* D_\mu \tilde{A_0}$ is positive.  Uniqueness 
of the right Polar Decomposition for the invertible matrix
$D_\lambda$ implies
$$D_\lambda = \tilde{A_0}^* D_\mu \tilde{A_0}.$$ Since the diagonal entries in
$D_\lambda$ and $D_\mu$ are listed in increasing order, it follows
that $D_\lambda=D_\mu$, hence $\rho_2$ is of the form $\rho_2(M) =
\sum_{i=1}^n \lambda_i^2 (g_i, M g_i)$.  Defining a unitary $U \in
M_n(\C)$ by letting $Ug_i = f_i$ for all $i$ and extending linearly,
we observe that
$$\rho_2(M) = \sum_{i=1}^n \lambda_i^2 (U^*f_i, M U^*f_i) =
\sum_{i=1}^n \lambda_i^2 (f_i, UMU^*f_i)= \rho_1(UMU^*)$$ for all $M
\in M_n(\C)$.  In other words, $\phi_2=(\phi_1)_U$.  \qed
\end{pf}

In \cite{bigpaper}, Powers constructed $E_0$-semigroups using
boundary weights over $L^2(0, \infty)$.  It is routine to check that
in our notation, these are the $E_0$-semigroups arising from the
boundary weight doubles $(\imath_\C, \eta)$, where $\imath_\C$ is
the identity map on $\C$ and $\eta$ is any boundary weight over
$L^2(0, \infty)$.
\begin{cor}\label{thesearenew}
Let $\phi: M_n(\C) \rightarrow M_n(\C)$ and $\psi: M_k(\C)
\rightarrow M_k(\C)$ be unital rank one $q$-pure maps, and let $\nu$
and $\eta$ be normalized unbounded boundary weights over $L^2(0,
\infty)$. Denote by $\alpha^d$ and $\beta^d$ the Bhat minimal
dilations of the $CP$-flows induced by the boundary weight doubles
$(\phi, \nu)$ and $(\psi, \mu)$, respectively.

If $n \neq k$, then $\alpha^d$ and $\beta^d$ are not cocycle
conjugate. In particular, if $n \neq 1$, then $\alpha^d$ is not
cocycle conjugate to the $E_0$-semigroup induced by $(\imath_\C,
\mu)$.
\end{cor}

\begin{pf}  From the proof of Theorem \ref{statesbig}, we know
that every corner $\gamma$ from $\phi$ to $\psi$ satisfies
$||\gamma||<1$ since $n \neq k$.  The result now follows from Lemma
\ref{generalnormone}. \qed \end{pf}

\section{Invertible unital $q$-pure maps}
Now that we have classified the unital $q$-pure maps on $M_n(\C)$ of
rank one, we explore the unital $q$-pure maps $\phi$ which are
invertible. In a stark contrast to the rank one case, we find that
for a given normalized unbounded boundary weight of the form
$\nu(\sqrt{I - \Lambda(1)} B \sqrt{I - \Lambda(1)})=(f, Bf)$ on
$L^2(0, \infty)$, the doubles $(\phi, \nu)$ and $(\psi, \nu)$
\textit{always} induce cocycle conjugate $E_0$-semigroups if $\phi$
and $\psi$ are unital invertible $q$-pure maps on $M_n(\C)$ and
$M_k(\C)$, respectively.
\newline \newline
The following proposition gives us a bijective correspondence
between invertible unital $q$-positive maps $\phi: M_n(\C)
\rightarrow M_n(\C)$ and unital conditionally negative maps $\Psi:
M_n(\C) \rightarrow M_n(\C)$:

\begin{prop}\label{cnegone}  If $\phi: M_n(\C) \rightarrow M_n(\C)$
is an invertible unital $q$-positive map, then $\phi^{-1}$ is
conditionally negative.  On the other hand, if $\Psi: M_n(\C)
\rightarrow M_n(\C)$ is a unital conditionally negative map, then
$\Psi$ is invertible and $\Psi^{-1}$ is $q$-positive.
\end{prop}
\begin{pf}  Let $\psi=\phi^{-1}$.  Since $\phi$ is self-adjoint, so
is $\psi$, and the first statement of the proposition now follows
from the fact that for large positive $t$ we have $$t \phi(I+t
\phi)^{-1} = t \psi^{-1}(I + t \psi^{-1})^{-1} = t (\psi+ tI)^{-1} =
\Big(I + \frac{\psi}{t}\Big)^{-1}=I- \frac{\psi}{t} +
\frac{\psi}{t}^2 - \ldots.$$

To prove the second statement, let $\Psi: M_n(\C) \rightarrow
M_n(\C)$ be any unital conditionally negative map.  Since $\Psi$ is
conditionally negative, it follows from a result of Evans and Lewis
in \cite{lewev} that $e^{-s \Psi}$ is completely positive for all $s
\geq 0$.  Therefore, $||e^{-s \Psi}||=||e^{-s\Psi}(I)||=
||e^{-s}I||= e^{-s}$ for all $s \geq 0$, and the integral $\int_0 ^
\infty e^{-s \Psi} ds$ converges. Observing that $(d/ds)
(-e^{-s \Psi}) = \Psi e^{-s \Psi}$, we find that
$$\Psi \Big(\int_0 ^ \infty e^{-s \Psi} ds\Big) = \int_0 ^ \infty
\Psi e^{-s \Psi}ds = \lim_{s \rightarrow \infty} (-e^{-s\Psi})\vert_0
^s = I,$$ so $\Psi$ is invertible and
$\Psi^{-1}= \int_0 ^ \infty
e^{-s \psi} ds$. Since $\Psi^{-1}$ is the integral of completely
positive maps, it is completely positive. 
Furthermore,  we find that $tI + \Psi$ is invertible for every $t>0$ and 
that $\Psi^{-1} 
\geq_q 0$, since the following holds for all $t>0$:
$$\int_0 ^\infty e^{-st}e^{-s \Psi}ds  = \int_0 ^\infty
e^{-s(tI+\Psi)}ds =
(tI + \Psi)^{-1} = \Psi^{-1}(I + t \Psi^{-1})^{-1}.$$ \qed 
\end{pf}

Examining the inverse of a unital invertible $q$-positive map $\phi$
is the key to finding the invertible $q$-subordinates of $\phi$, as
we find in the following proposition and corollary:

\begin{prop}\label{invsubs}  Let $\phi_1: M_n(\C) \rightarrow M_n(\C)$
be an invertible unital $q$-positive map, and let
$\psi_1=\phi_1^{-1}$. Suppose $\psi_2: M_n(\C) \rightarrow M_n(\C)$
is conditionally negative and $\psi_2-\psi_1$ is completely
positive. Then $\psi_2$ is invertible, and $\phi_2:=(\psi_2)^{-1}$
satisfies $\phi_1 \geq_q \phi_2 \geq_q 0$.
\end{prop}
\begin{pf} Assume the hypotheses of the proposition, and let $s>0$
be arbitrary.  Define a function $f$ on $\R$ by $f(t) = e^{-ts
\psi_1}e^{(t-1)s\psi_2}$.  The equality below is $f(1)-f(0)=\int_0^1
f'(t) dt$:
$$e^{-s\psi_1}-e^{-s\psi_2} =  \int_0^1 s e^{-ts\psi_1}(\psi_2-\psi_1)
e^{(t-1)s\psi_2} dt.$$ The inside of the integral above
is the composition of completely positive maps, so
$e^{-s\psi_1}-e^{-s\psi_2}$ is completely positive.  This implies
$e^{-s \psi_1}(I) - e^{-s \psi_2}(I) \geq 0$, so
$$||e^{-s \psi_2}||=||e^{-s \psi_2}(I)|| \leq ||e^{-s \psi_1}(I)|| =
||e^{-s}(I)||=e^{-s}.$$  Now the argument given in the previous
proposition shows that $\int_0 ^\infty e^{-s \psi_2} ds$ converges
and is equal to $\psi_2^{-1}$.  Letting $\phi_2 = \psi_2^{-1}$, we
observe that $\phi_1 \geq_q \phi_2$ since the quantity below is 
completely positive for every $t \geq 0$:
$$\phi_1(I + t \phi_1)^{-1} - \phi_2(I + t \phi_2)^{-1}
= \int_0 ^\infty e^{-st}(e^{-s\psi_1} - e^{-s \psi_2}) ds.$$ \qed \end{pf}

\begin{cor}\label{invsubsone} Let $\phi_1: M_n(\C) \rightarrow M_n(\C)$
be an invertible unital $q$-positive map, and let $\phi_2: M_n(\C)
\rightarrow M_n(\C)$ be linear and invertible.

Then $\phi_1 \geq_q \phi_2 \geq_q 0$ if and only if $\phi_2^{-1}$ is
conditionally negative and $\phi_2^{-1} - \phi_1^{-1}$ is completely
positive.
\end{cor}

\begin{pf}  The backward direction follows from Proposition
\ref{invsubs}. Assume the hypotheses of the forward direction and
let $\psi_1=\phi_1^{-1}$ and $\psi_2= \phi_2^{-1}$.  Since
$\phi_2$ is self-adjoint, so is $\psi_2$.  For sufficiently large
positive $t$ we have
$$t \phi_2(I + t \phi_2)^{-1} = \Big(I + \frac{\psi_2}{t}\Big)^{-1}=I-
\frac{\psi_2}{t} + \frac{\psi_2^2}{t^2} - \ldots$$ and
$$t^2(\phi_1(I + t \phi_1)^{-1} - \phi_2(I + t \phi_2)^{-1}) = \psi_2 -
\psi_1 + \Big(\frac{\psi_2^2-\psi_1^2}{t} - \frac{\psi_2^3 -
\psi_1^3}{t^2} + \ldots\Big).$$  The first equation shows that
$\phi_2^{-1}$ is conditionally negative, while the second shows that
$\phi_2^{-1} - \phi_1^{-1}$ is completely positive.  \qed
\end{pf}

Now that we know how to find all invertible $q$-subordinates of an
invertible unital $q$-positive map $\phi$, we ask if there can be
any other $q$-subordinates of $\phi$. We will find that the answer
is no (see Proposition \ref{invert}).  Proving this will require the
use of some machinery (notably Lemma \ref{eps}), which we now build.

\begin{defn}\label{defeps}  For every $\phi: M_n(\C) \rightarrow M_n(\C)$ and
$ \epsilon \in [0,1]$, we define a map $\phi_\epsilon$ by
$\phi_\epsilon = \epsilon I + (1- \epsilon)\phi$.
\end{defn} If $\phi$ is $q$-positive, then $\phi_\epsilon$ is
invertible for all $\epsilon \in (0,1]$.  In the lemmas that follow,
we make frequent use of the fact that for all $t \geq 0$ we have
\begin{eqnarray}\label{I-} t \phi(I + t \phi)^{-1} = I - (I + t
\phi)^{-1}.\end{eqnarray}  We present a quick consequence of
\eqref{I-} for all $a \geq 0$ and $b \geq 0$:
\begin{eqnarray} \label{I-better} a(I + bt \phi)^{-1}=aI - abt\phi(I + bt \phi)^{-1}
\end{eqnarray}

\begin{lem}  Let $\phi: M_n(\C) \rightarrow M_n(\C)$ be completely positive.
If $\phi_{\epsilon_k} \geq_q 0$ for some monotonically decreasing
sequence $\{\epsilon_k\}$ of positive real numbers tending to $0$,
then $\phi \geq_q 0$.
\end{lem}
\begin{pf}  Assume the hypotheses of the lemma.  Let $k$ be
arbitrary.  Since $\phi_{\epsilon_k} \geq_q 0$, we know $I- (I
+ t \phi_{\epsilon_k})^{-1}$ is completely positive for all $t \geq
0$. Noting that
$$I - (I + t \phi_\epsilon)^{-1}  = I - \Big((1 +
t \epsilon I) + (1-\epsilon)t \phi \Big)^{-1}= I - \frac{1}{1+t
\epsilon} \Big(I + \frac{t(1- \epsilon)}{1+t \epsilon} \phi
\Big)^{-1}$$ and substituting $t'=t(1-\epsilon_k)/(1+t
\epsilon_k)$, we see
$$I - (I + t \phi_\epsilon)^{-1}=
I - \frac{1}{1+ (\frac{\epsilon_k}{1-\epsilon_k + t'\epsilon_k})t'}
(I + t' \phi)^{-1}.$$  Varying $t$ throughout $[0, \infty)$, we find
that the above equation is completely positive for all $t' \in [0,
-1+1/\epsilon_k)$.  Of course, for any $t' \in [0,
-1+1/ \epsilon_k )$, we have $t' \in [0, -1+1/\epsilon_\ell)
$ for all $ \ell \geq k$ by
monotonicity of the sequence $\{\epsilon_n\}$.  Therefore, we may
repeat the same argument to conclude that for any $t' \in [0,
-1+1/\epsilon_k)$, the map $$I - \frac{1}{1+
(\frac{\epsilon_\ell}{1-\epsilon_\ell + t'\epsilon_\ell})t'} (I + t'
\phi)^{-1}$$ is completely positive for all $\ell \geq k$.

Now fix any $t' > 0$, so $t' \in (0, -1+1/\epsilon_k)$ for some $k \in \mathbb{N}$.  A
straightforward computation shows that the sequence $\{c_n\}$
defined by $c_n = \epsilon_n / (1-\epsilon_n + t' \epsilon_n)$
monotonically decreases to $0$. From the previous paragraph, we know
that the map
$$I - \frac{1}{1+ c_\ell t'} (I + t'
\phi)^{-1}$$ is completely positive for all $\ell \geq k$.  Since
$c_n \downarrow 0$ it follows that
$$I - (I + t' \phi)^{-1}$$ is completely positive.  In other words,
$t' \phi(I + t' \phi)^{-1}$ is completely positive.  Since $t'>0$ was
chosen arbitrarily and $\phi$ is completely positive, the lemma
follows. \qed
\end{pf}

\begin{lem} If $\phi: M_n(\C) \rightarrow M_n(\C)$
and $\phi \geq_q 0$, then $\phi_\epsilon \geq _q 0$ for all
$\epsilon \in [0,1)$.
\end{lem}

\begin{pf}  Suppose that $\phi \geq_q 0$, and let $\epsilon \in
[0,1)$ be arbitrary.  For each $t>0$, we apply formula
\eqref{I-better} to $a=1/(1+t\epsilon)$ and
$b=t(1-\epsilon)/(1+t \epsilon)$ to find
\begin{eqnarray*}I - (I + t \phi_\epsilon)^{-1} & = & I - \frac{1}{1+t
\epsilon} \Big(I + \frac{t(1- \epsilon)}{1+t \epsilon} \phi
\Big)^{-1}  \\ & = & \Big(1-\frac{1}{1+t\epsilon}\Big)I +
\frac{t(1-\epsilon)}{(1+t \epsilon)^2} \phi \Big(I +
\frac{t(1-\epsilon)}{1+t \epsilon} \phi\Big)^{-1},
\end{eqnarray*} where both terms on the last line are completely
positive by assumption.  Furthermore, $\phi_\epsilon$ is completely
positive, hence $\phi_\epsilon \geq_q 0$. \qed \end{pf}

\begin{cor}  Let $\phi: M_n(\C) \rightarrow M_n(\C)$ be a completely positive map.
Then $\phi \geq_q 0$ if and only if $\phi_\epsilon \geq_q 0$ for all
$\epsilon \in (0,1)$.
\end{cor}

\begin{lem}\label{eps}  Let $\phi: M_n(\C) \rightarrow M_n(\C)$ and $\psi: M_n(\C)
\rightarrow M_n(\C)$ be $q$-positive maps.  Then $\phi \geq_q \psi$
if and only if $\phi_\epsilon \geq_q \psi_\epsilon$ for all
$\epsilon \in (0,1)$.
\end{lem}

\begin{pf}  For any $\epsilon \in (0,1)$ we have $\phi_\epsilon -
\psi_\epsilon = \epsilon(\phi-\psi)$, so $\phi- \psi$ is completely
positive if and only if $\phi_\epsilon - \psi_\epsilon$ is
completely positive for all $\epsilon \in (0,1)$.  For all $t'>0$ we
have \begin{eqnarray}\label{short} t' \Big( \phi(I + t' \phi)^{-1} -
\psi(I + t' \psi)^{-1} \Big) = (I + t' \psi)^{-1} - (I + t'
\phi)^{-1},\end{eqnarray} and for all $t>0$ we have

\begin{eqnarray} t (\phi_\epsilon(I + t \phi_\epsilon)^{-1}  -  \psi_\epsilon(I +
t \psi_\epsilon)^{-1})  = \Big(I - (I + t \phi_\epsilon)^{-1}\Big) -
\Big(I - (I + t \psi_\epsilon)^{-1} \Big) \nonumber \\ \label{long}
 =  \frac{1}{1+t \epsilon} \Big((I + \frac{t(1-\epsilon)}{1+t
\epsilon} \psi)^{-1} - (I + \frac{t(1-\epsilon)}{1+t \epsilon}
\phi)^{-1}\Big).\end{eqnarray}

Assume the hypotheses of the forward direction.  Showing that
$\phi_\epsilon \geq_q \psi_\epsilon$ for all $\epsilon \in (0,1)$ is
equivalent to proving that \eqref{long} is completely positive for
every $t \in (0, \infty)$ and $\epsilon \in (0,1)$.  But this
follows from complete positivity of $\eqref{short}$ since
$t(1-\epsilon)/(1+t \epsilon) \in (0, \infty)$ for every
$\epsilon \in (0, 1)$ and $t \in (0, \infty)$.  Now assume the
hypotheses of the backward direction.  Any $t' \in (0, \infty)$ can
be written as $t(1-\epsilon)/(1+t \epsilon)$ for some $\epsilon
\in (0,1)$ and $t \in (0, \infty)$, so complete positivity of
\eqref{long} for all such $\epsilon$ and $t$ implies that
\eqref{short} is completely positive for all $t' >0$, hence $\phi
\geq_q \psi$. \qed \end{pf}

We are now in a position to prove what is perhaps the most striking
result of the section:
\begin{prop}\label{invert}
Let $\xi: M_n(\C) \rightarrow M_n(\C)$ be an invertible unital
$q$-positive map. If $\phi: M_n(\C) \rightarrow M_n(\C)$ is
$q$-positive and $\xi \geq_q \phi$, then $\phi$ is either invertible
or identically zero.
\end{prop}

\begin{pf}  For every $\epsilon \in (0,1)$, form $\xi_\epsilon$ and
$\phi_\epsilon$ as in Definition \ref{defeps}, and let
$\psi_\epsilon: = (\phi_\epsilon)^{-1}$.  By Lemma \ref{eps}
we have $\xi_{\epsilon} \geq_q \phi_{\epsilon}$ for each $\epsilon$, so 
$\psi_\epsilon$
is conditionally negative and $\psi_\epsilon - (\xi_\epsilon)^{-1}$ 
is completely positive by Corollary \ref{invsubsone}. We first
examine the case when the norms $||\psi_\epsilon||$ remain bounded
as $\epsilon \rightarrow 0$. More precisely, suppose that for all
$\epsilon$ sufficiently small we have $||\psi_\epsilon|| < r$ for
some $r >0$.  By compactness of the closed unit ball of radius $r$
in $B(M_n(\C))$, there is a decreasing sequence $\{\epsilon_k\}_{k
\in \mathbb{N}}$ converging to $0$ such that
$\{\psi_{\epsilon_k}\}_{k \in \mathbb{N}}$ has a (bounded) norm
limit $\psi$ as $k \rightarrow \infty$.  Noting that
$$I- \phi\psi = \phi_{\epsilon_k}\psi_{\epsilon_k} - \phi \psi =
(\phi_{\epsilon_k} - \phi)(\psi_{\epsilon_k}-\psi) +
\phi(\psi_{\epsilon_k} - \psi) + (\phi_{\epsilon_k} - \phi)\psi$$
and then applying the triangle inequality, we find that
\begin{eqnarray*}
||I - \phi \psi|| &=& ||\phi_{\epsilon_k}\psi_{\epsilon_k} - \phi
\psi|| \\ & \leq & ||\phi_{\epsilon_k} - \phi|| \
||\psi_{\epsilon_k} -\psi|| + ||\phi|| \ ||\psi_{\epsilon_k}-\psi||
+ ||\phi_{\epsilon_k}-\phi|| \ ||\psi||
\end{eqnarray*} for all $k \in \mathbb{N}$.  But $\phi$ and $\psi$ are
bounded maps while $\psi_{\epsilon_k}\rightarrow \psi$ in norm and
$\phi_{\epsilon_k} \rightarrow \phi$ in norm, so the above equation
tends to $0$ as $k \rightarrow \infty$. We conclude that $\phi \psi
= I$.  Similarly $\psi \phi = I$, hence $\phi$ is invertible and
$\psi=\phi^{-1}$.

If the first case does not hold, then for some decreasing sequence
$\{\epsilon_k\}$ tending to zero, the norms
$\{||\psi_{\epsilon_k}||\}_{k \in \mathbb{N}}$ form an unbounded
sequence. For each $k \in \mathbb{N}$, we write
$$(\xi_{\epsilon_k})^{-1}(A)= s_k A + Y_k A + AY_k^* -
\sum_{i=1}^{m_k} S_{k_i} A S_{k_i}^*$$ $$\textrm{ and }$$
$$\psi_{\epsilon_k}(A)= t_k A + Z_k A + AZ_k^* -
\sum_{i=1}^{\ell_k}T_{k_i} A T_{k_i}^*,$$ where $m_k, \ell_k \leq n^2$,
$s_k \in \R$, $t_k
\in \R$, $tr(Y_k)=tr(Z_k)=0$, 
$tr(S_{k_i})=0$ and $tr(S_{k_i}^*S_{k_j})$ is non-zero if and only if $i=j$ ($i,j \leq m_k)$,
and $tr(T_{k_i})=0$ and $tr(T_{k_i}^* T_{k_j})$ is non-zero if and only if $i=j$ ($i,j \leq \ell_k$).

Since $\psi_{\epsilon_k} -
(\xi_{\epsilon_k})^{-1}$ is completely positive for all $k \in
\mathbb{N}$, we know that for each $k$,
there exist $p_k \leq n^2$, complex numbers $\{x_{k_i}\}_{i=1}^{p_k}$, 
and maps $\{X_{k_i}\}_{i=1}^{p_k}$ with $tr(X_{k_i})=0$, such
that for all $A \in M_n(\C)$,

\begin{eqnarray} (\psi_{\epsilon_k}  - (\xi_{\epsilon_k})^{-1})(A) &
= &
\sum_{i=1}^{p_k} (X_{k_i}+x_{k_i} I) A (X_{k_i} + x_{k_i}I)^* \nonumber \\
& = &  \Big(\sum_{i=1}^{p_k} |x_{k_i}|^2\Big)A +
\Big(\sum_{i=1}^{p_k} \overline{x_{k_i}}X_{k_i}\Big)A+
A\Big(\sum_{i=1}^{p_k} \overline{x_{k_i}}X_{k_i}\Big)^* \nonumber \\
\label{l1} & \ & \ + \sum_{i=1}^{p_k} X_{k_i}AX_{k_i}^*.
\end{eqnarray}
Simultaneously, for all $A \in M_n(\C)$ we have
\begin{eqnarray} (\psi_{\epsilon_k} - (\xi_{\epsilon_k})^{-1})(A) &
= & (t_k -s_k)A + (Z_k-Y_k)A +A(Z_k-Y_k)^* \nonumber \\ \label{l2} &
\ & + \Big(\sum_{i=1}^{m_k} S_{k_i} A S_{k_i}^* -
\sum_{i=1}^{\ell_k}T_{k_i} A T_{k_i}^* \Big).
\end{eqnarray}
We claim that \begin{eqnarray}\label{Xis} \Big|\Big|\sum_{i=1}^{p_k}
X_{k_i}AX_{k_i}^*\Big|\Big| \leq \Big|\Big| \sum_{i=1}^{m_k} S_{k_i}
A S_{k_i}^*\Big|\Big|\end{eqnarray} for all $k \in \mathbb{N}$.  To
prove this, we let $\{v_j\}_{j=1}^n$ be any orthonormal basis for
$\C^n$, let $h_j=v_j / \sqrt{n}$ for each $i$, let $f \in \C^n$ be
arbitrary, and define maps $A_j$ for $j=1, \ldots, n$ by $A_j = f
h_j^*$.  Using the trace conditions on the maps $Y_k$, $Z_k$,
$\{T_{k_i}\}$, $\{S_{k_i}\}$, and $\{X_{k_i}\}$, we find that

\begin{eqnarray*} \sum_{j=1}^n(\psi_{\epsilon_k} - (\xi_{\epsilon_k})^{-1})(A_j)
h_j & = & (t_k -s_k) f + (Z_k - Y_k)f \\ & = & \Big(\sum_{i=1}^{p_k}
|x_{k_i}|^2\Big)f + \Big(\sum_{i=1}^{p_k} \overline{x_{k_i}}
X_{k_i}\Big)f.\end{eqnarray*} Since $f$ was arbitrary, it follows
that $$\Big(t_k - s_k - \sum_{i=1}^{p_k} |x_{k_i}|^2\Big)I =
\Big(\sum_{i=1}^{p_k} \overline{x_{k_i}} X_{k_i} \Big) -
(Z_k-Y_k).$$ Taking the trace of both sides yields
$$0 = tr\Big(\Big(\sum_{i=1}^{p_k} \overline{x_{k_i}} X_{k_i} \Big) - (Z_k-Y_k)\Big) =
tr\Big((t_k - s_k -\sum_{i=1}^{p_k} |x_{k_i}|^2)I\Big),$$ so $t_k -
s_k = \sum_{i=1}^{p_k} |x_{k_i}|^2$ and $Z_k - Y_k =
\sum_{i=1}^{p_k} \overline{x_{k_i}}X_{k_i}$.  Formulas \eqref{l1}
and \eqref{l2} now imply that $$\sum_{i=1}^{p_k} X_{k_i}AX_{k_i}^*
=\Big(\sum_{i=1}^{m_k} S_{k_i} A S_{k_i}^* -
\sum_{i=1}^{\ell_k}T_{k_i} A T_{k_i}^* \Big).$$  Therefore, the map
$A \rightarrow \sum_{i=1}^{m_k} S_{k_i} A S_{k_i}^* -
\sum_{i=1}^{\ell_k}T_{k_i} A T_{k_i}^*$ is completely positive, and
$$\Big|\Big|\sum_{i=1}^{p_k} X_{k_i}X_{k_i}^*\Big|\Big| =
\Big|\Big|\sum_{i=1}^{m_k} S_{k_i}S_{k_i}^* -
\sum_{i=1}^{\ell_k}T_{k_i}T_{k_i}^*\Big|\Big| \leq
\Big|\Big|\sum_{i=1}^{m_k} S_{k_i}S_{k_i}^*\Big|\Big|,$$
establishing
\eqref{Xis}.

We now show that there exists some $M \in \mathbb{N}$ such that
\begin{eqnarray}\label{Xibound}||X_{k_i}|| \leq M\end{eqnarray}
for all $k \in \mathbb{N}$ and $i \in \{1, \ldots, p_k\}$.  To do
this, we first note that since the sequence of invertible maps
$\{\xi_{\epsilon_k}\}_{k \in \mathbb{N}}$ converges in norm to the
invertible map $\xi$, the sequence $\{(\xi_{\epsilon_k})^{-1}\}_{k
\in \mathbb{N}}$ converges in norm to $\xi^{-1}$.  Write $\xi^{-1}$
in the form
$$\xi^{-1}(A) = sA + YA + AY^* - \sum_{i=1}^m S_i A S_i^*,$$ where $m \leq n^2$,
$s \in \R$,
$tr(Y)=0$, and for all $i$ and $j$, $tr(S_i)=0$ and $tr(S_iS_j^*)$ is non-zero if 
and only if $i=j$.
Let $f \in \C^n$ be arbitrary, and define vectors $\{h_j\}_{j=1}^n$
and maps $\{A_j\}_{j=1}^n$ exactly as we did earlier in the proof.  Then
$\sum_{j=1}^n (\xi_{\epsilon_k})^{-1}(A_j)h_j = s_k f + Y_kf$ for all $k \in \mathbb{N}$
and $\sum_{j=1}^n \xi^{-1}(A_j)h_j = s f + Yf$.  Since $(\xi_{\epsilon_k})^{-1}$
converges to $\xi^{-1}$ as $k \rightarrow \infty$, we see that
$(s_k - s) f + (Y_k - Y)f$ 
converges to $0$ as $k \rightarrow \infty$.  But $f$ was arbitrary, so
$$\lim_{k \rightarrow \infty} \Big((s_k - s) I + Y_k - Y \Big)= 0.$$  The 
limit of the trace
of the above equation must also be zero, so $s_k$ converges to $s$
and consequently $Y_k$ converges to $Y$.
This implies that not only are the sequences of complex numbers $\{s_k\}_{k=1}^\infty$
and maps $\{Y_k\}_{k=1}^\infty$ both bounded, but that 
the sequence of linear maps $\{W_k\}_{k=1}^\infty$ defined by
$W_k(A)= \sum_{i=1}^{m_k} S_{k_i}AS_{k_i}^*$
is bounded and converges to the map $W(A) = \sum_{i=1}^m S_i A S_i^*$. 
Choose $M \in \mathbb{N}$ so that $M^2 \geq n^2 \sup_{k
\in \mathbb{N}} \{||W_k||\}$.  For every $k \in \mathbb{N}$ and $i
\in \{1, \ldots, m_k\}$, we have $||S_{k_i}||^2 \leq ||W_k|| \leq
M^2/n^2$.  Combining this fact with \eqref{Xis}, we find that for
every $k \in \mathbb{N}$ and $i \in \{1, \ldots, p_k\}$,
\begin{eqnarray*}||X_{k_i}||^2 & = & ||X_{k_i}X_{k_i}^*|| \leq
||\sum_{i=1}^{p_k} X_{k_i}X_{k_i}^*|| \leq ||\sum_{i=1}^{m_k}
S_{k_i}S_{k_i}^*|| \leq  \sum_{i=1}^{m_k}||S_{k_i}||^2
\\ & \leq & n^2 \max\{||S_{k_i}||^2: i=1, \dots, m_k\} \leq M^2,
\end{eqnarray*} proving \eqref{Xibound}.

Since $||\psi_{\epsilon_k}||
\rightarrow \infty$ as $k \rightarrow \infty$ while
$||(\xi_{\epsilon_k})^{-1}|| \rightarrow ||(\xi)^{-1}|| < \infty$,
there is a sequence of maps $\{A_{\epsilon_k}\}$ of norm one such
that
$||(\psi_{\epsilon_k}-(\xi_{\epsilon_k})^{-1})(A_{\epsilon_k})||
\rightarrow \infty$ as $k \rightarrow \infty$.  However, we also
have
\begin{eqnarray}||(\psi_{\epsilon_k}-(\xi_{\epsilon_k})^{-1})(A_{\epsilon_k})||
& = &\Big|\Big| \Big(\sum_{i=1}^{p_k} |x_{k_i}|^2\Big)A_{\epsilon_k}
+ \Big(\sum_{i=1}^{p_k}
\overline{x_{k_i}}X_{k_i}\Big)A_{\epsilon_k} \nonumber \\ & \ & \ +
A_{\epsilon_k}\Big(\sum_{i=1}^{p_k} \overline{x_{k_i}}X_{k_i}\Big)^*
+ \sum_{i=1}^{p_k} X_{k_i}A_{\epsilon_k}X_{k_i}^*\Big|
\Big| \nonumber \\ 
& \label{above} \leq & \sum_{i=1}^{p_k} |x_{k_i}|^2 + 2 M \sum_{i=1}^{p_k}
|x_{k_i}| + p_k M^2.
\end{eqnarray}
We note that
\begin{eqnarray} \label{prob} \Big(\sum_{i=1}^{p_k}|x_{k_i}|\Big)^2 \geq
\sum_{i=1}^{p_k}|x_{k_i}|^2 \geq
\frac{(\sum_{i=1}^{p_k}|x_{k_i}|)^2}{p_k} \geq
\frac{(\sum_{i=1}^{p_k}|x_{k_i}|)^2}{n^2}\end{eqnarray} for all $k$.
For each $k$, let $\lambda_k = \sum_{i=1}^{p_k} |x_{k_i}|$, noting that
$\lambda_k \rightarrow \infty$ as $k \rightarrow
\infty$ since Eq.\eqref{above} tends to infinity as $k \rightarrow \infty$.
 Let $A
\in M_n(\C)$ be any matrix such that $||A||=1$, and let $C= \sup_{k
\in \mathbb{N}} ||(\xi_{\epsilon_k})^{-1}||< \infty$. Using the
reverse triangle inequality and \eqref{prob}, we find that for each
$k \in \mathbb{N}$,
\begin{eqnarray}||\psi_{\epsilon_k}(A)|| & \geq &
||(\psi_{\epsilon_k}-(\xi_{\epsilon_k})^{-1})(A)|| -
||(\xi_{\epsilon_k})^{-1}(A)|| \nonumber \\ \label{last} & \geq &
\frac{\lambda_k^2}{n^2} - 2M \lambda_k - n^2 M^2-C.\end{eqnarray}
Since $\lim_{k \rightarrow \infty} \lambda_k = \infty$, 
Eq.\eqref{last} tends to infinity as $k \rightarrow \infty$.  For all
$k$ large enough that Eq.\eqref{last} is positive, we have
$$||\phi_{\epsilon_k}|| =\frac{1}{\inf \{||\psi_{\epsilon_k}(A)||: ||A||=1
\}} \leq \frac{1}{\lambda_k^2 / n^2 - 2M \lambda_k - n^2 M^2-C},$$
so $\lim_{k \rightarrow \infty} ||\phi_{\epsilon_k}|| = 0$. But the
sequence $\{\phi_{\epsilon_k}\}_{k=1}^\infty$ converges to $\phi$ in
norm, hence $\phi \equiv 0$.  \qed \end{pf}

\begin{prop}\label{cnegpure} An invertible unital linear map $\phi: M_n(\C)
\rightarrow M_n(\C)$ is $q$-pure if and only if $\phi^{-1}$ is of
the form
$$\phi^{-1}(A) = A + YA + AY^*$$ for some $Y=-Y^* \in M_n(\C)$ such
that $tr(Y)=0$.
\end{prop}

\begin{pf} Let $\psi=\phi^{-1}$.  Assume the hypotheses of the
forward direction. Write
$$\psi(A) = sA + YA + AY^* - \sum_{i=1}^k \lambda_i X_i A X_i^*,$$
where $s \in \R$, $tr(Y)=0$, and for each $i$ and $j$ we have
$\lambda_i \geq 0$, $tr(X_i)=0$, and $tr(X_i^*X_j)= n \delta_{ij}$.

Defining $\psi': M_n(\C) \rightarrow M_n(\C)$ by
$$\psi'(A) = sA+YA+AY^*,$$
we note that $\psi'$ is conditionally negative, and $\psi'-\psi$ is
completely positive since $(\psi'-\psi)(A)= \sum_{j=1}^k \lambda_j
X_j A X_j^*$ for all $A$.  By Lemma \ref{invsubs}, it follows that
$\psi'$ is invertible and that $\phi':=(\psi')^{-1}$ satisfies $\phi
\geq_q \phi' \geq_q 0$.

Since $\phi$ is $q$-pure, there is some $t_0 \geq 0$ such that
$\phi' = \phi^{(t_0)}$, hence
$$\psi'=(\phi')^{-1}=\Big(\phi(I + t_0 \phi)^{-1}\Big)^{-1} =\Big(\psi^{-1}(I + t_0 \psi^{-1})\Big)^{-1}= \Big((t_0
I + \psi)^{-1}\Big)^{-1}= t_0I + \psi.$$  Therefore, for all $A \in
M_n(\C)$ we have
$$\psi'(A) =\psi(A)+ \sum_{j=1}^k
\lambda_j X_j A X_j^* = \psi(A) + t_0A,$$ so the map $L: A
\rightarrow \lambda_j X_j A X_j^*$ satisfies $L=t_0 I$.  We repeat a
familiar argument:  Let $f \in \C^n$ be arbitrary, choose an
orthonormal basis $\{v_k\}_{k=1}^n$ of $\C^n$, define $h_k = v_k/
\sqrt{n}$ for each $k$, and form $\{A_k\}_{k=1}^n$ by $A_k =
fh_k^*$. The trace conditions for the maps $\{X_j\}$ imply that
$\sum_{k=1}^n L(A_k)h_k=0$.  However, since $L=t_0 I$, we must also
have $\sum_{k=1}^n L(A_k)h_k = t_0 f$. From arbitrariness of $f$, we
conclude $t_0=0$. Therefore, $\psi$ has the form $\psi(A)=sA + YA +
AY^*$. Since $\psi(I)=I=sI+Y+Y^*$ and $tr(Y)=0$, we have $s=1$ and
consequently $Y=-Y^*$.

Now assume the hypotheses of the backward direction.  Note that
$\psi$ is conditionally negative and unital, hence $\phi$ is
$q$-positive by Proposition \ref{cnegone}. Let $\Phi$ be any non-zero
$q$-positive map such that $\phi \geq_q \Phi$, so by Corollary
\ref{invsubsone} and Proposition \ref{invert}, $\Phi$ is invertible
and $\Psi:=(\Phi)^{-1}$ is a conditionally negative map such that
$\Psi-\psi$ is completely positive. Write $\Psi$ in the form
$$\Psi(A) = s'A+ ZA + AZ^* - \sum_{i=1}^{m} \mu_i T_i A T_i^*,$$ where
$s' \in \R$ and for all $i$ and $j$, $\mu_i > 0$, $tr(T_i)=0$, and
$tr(T_i^*T_j)=n \delta_{ij}$.  Writing $C=Z-Y$, we have
$$(\Psi-\psi)(A)= (s'-1)A + CA+AC^* - \sum_{i=1}^m \mu_i T_i A T_i^*.$$
By a familiar argument, complete positivity of $\Psi - \psi$ and the trace conditions for
the above maps imply that $s' \geq 1$, $C=0$, and $T_i=0$ for all
$i$. Therefore $\Psi = \psi + (s'-1)I$, so $\Phi=\Psi^{-1}=
\phi^{(s'-1)}$. We conclude that $\phi$ is $q$-pure. \qed \end{pf}

Let the matrices $\{e_{jk}\}_{j,k=1}^n$ denote the standard basis for $M_n(\C)$,
writing each $A= (a_{jk}) \in M_n(\C)$ as $A= \sum_{j,k} a_{jk}e_{jk}$.  The
following theorem classifies all unital invertible $q$-pure maps on $M_n(\C)$:

\begin{thm}\label{phiu}An invertible unital linear map $\phi: M_n(\C)
\rightarrow M_n(\C)$ is $q$-pure if and only if for some unitary $U
\in M_n(\C)$, the map $\phi_U$ is the Schur map

\begin{equation*} \phi_U(a_{jk}e_{jk}) = \left\{
\begin{array}{cc}
\frac{a_{jk}}{1+i(\lambda_j - \lambda_k)}e_{jk} & \textrm{if } j<k \\
a_{jk}e_{jk} & \textrm{if } j=k \\
\frac{a_{jk}}{1-i(\lambda_j - \lambda_k)}e_{jk}& \textrm{if } j>k
\end{array} \right.
\end{equation*}
for all $A= (a_{jk}) \in M_n(\C)$ and $j,k=1, \ldots, n$, where $\lambda_1, \ldots, \lambda_n \in
\R$ and $\lambda_1 + \ldots + \lambda_n = 0$.
\end{thm}

\begin{pf} Assume the hypotheses of the forward direction.  By the previous
proposition, $\psi:=\phi^{-1}$ has the form $\psi(A) = A +
\tilde{Y}A + A \tilde{Y}^*$ for some $\tilde{Y} \in M_n(\C)$ with
$\tilde{Y}=-\tilde{Y^*}$ and $tr(\tilde{Y})=0$. Let $B=
-i\tilde{Y}$, so $B=B^*$.  Defining $Y:= (1/2)I
+\tilde{Y}=(1/2)I+iB$, we find $\psi(A)= YA+AY^*$ for all
$A\in M_n(\C)$.  Since $B$ is self-adjoint, there is some unitary $U
\in M_n(\C)$ such that $U^* B U$ is a diagonal matrix $D$.  For each
$k \in \{1, \ldots, n\}$ let $\lambda_k \in \R$ be the $kk$ entry of
$D$.  Note that since $tr(B)=0$ we have $\sum_{k=1}^n \lambda_k =
0$, and that $U^*YU$ is the diagonal matrix $M$ whose $kk$ entry is
$1/2 +i \lambda_k$. Defining a map $\psi_U$ by
$\psi_U(A)=U^*\psi(UAU^*)U$ for all $A \in M_n(\C)$, we find that

\begin{eqnarray*}\psi_U(A) &=& U^*(YUAU^* + UAU^*Y^*)U
\\ &=&(U^*YU)A+A(U^*YU)^* = MA+AM^*.\end{eqnarray*}
A quick calculation shows that this is just the Schur map
\begin{equation*} \psi_U(a_{jk}e_{jk}) = \left\{
\begin{array}{cc}
({1+i(\lambda_j - \lambda_k)})a_{jk}e_{jk} & \textrm{if } j<k \\
a_{jk}e_{jk} & \textrm{if } j=k \\
({1-i(\lambda_j - \lambda_k)})a_{jk}e_{jk}& \textrm{if } j>k
\end{array} \right.,
\end{equation*}
and so $(\psi_U)^{-1}$ has the form
\begin{equation*} (\psi_U)^{-1}(a_{jk}e_{jk}) = \left\{
\begin{array}{cc}
\frac{a_{jk}}{1+i(\lambda_j - \lambda_k)}e_{jk} & \textrm{if } j<k \\
a_{jk}e_{jk} & \textrm{if } j=k \\
\frac{a_{jk}}{1-i(\lambda_j - \lambda_k)}e_{jk}& \textrm{if } j>k
\end{array} \right.. \end{equation*} It is straightforward to verify
that $(\psi_U)^{-1}$ is the map $\phi_U(A)= U^*\phi(UAU^*)U$.

Assume the hypotheses of the backward direction. Let $T$ be the
diagonal matrix whose $kk$th entry is $\lambda_k$ for every $k=1,
\ldots, n$.  We observe that $tr(T)=0$ and $T=T^*$.  Now let $C=iT$,
and let $\tilde{T}= (1/2)I +C$.  We routinely verify that
$C=-C^*$ and $tr(C)=0$, and that $(\phi_U)^{-1}$ satisfies
$(\phi_U)^{-1}(A)=\tilde{T}A+A\tilde{T}^* = A + CA + AC^*$ for all
$A \in M_n(\C)$.  Proposition \ref{cnegpure} implies that $\phi_U$
is $q$-pure, whereby $\phi$ is $q$-pure by Proposition
\ref{basischange}.
 \qed \end{pf}

As it turns out, boundary weight doubles $(\phi, \nu)$ for invertible unital $q$-pure maps
$\phi: M_n(\C) \rightarrow M_n(\C)$ and normalized unbounded boundary weights
$\nu$ over $L^2(0,
\infty)$ of the form $\nu(\sqrt{I - \Lambda(1)} B \sqrt{I -
\Lambda(1)}) = (f,Bf)$ give us
nothing new in terms of $E_0$-semigroups:

\begin{thm}  Let $\phi: M_n(\C) \rightarrow M_n(\C)$ be unital, invertible, and
$q$-pure, and let $\nu$ be a normalized unbounded boundary weight
over $L^2(0,\infty)$ of the form $\nu(\sqrt{I - \Lambda(1)} B
\sqrt{I - \Lambda(1)}) = (f,Bf)$. Then $(\phi, \nu)$ and
$(\imath_\C, \nu)$ induce cocycle conjugate $E_0$-semigroups.
\end{thm}

\begin{pf}  By Theorem \ref{phiu} and Propositions
\ref{basischange} and \ref{hypqc}, we may assume that $\phi$ is the
Schur map
\begin{equation*} \phi(a_{jk}e_{jk}) = \left\{
\begin{array}{cc}
\frac{a_{jk}}{1+i(\lambda_j - \lambda_k)}e_{jk} & \textrm{if } j<k \\
a_{jk}e_{jk} & \textrm{if } j=k \\
\frac{a_{jk}}{1-i(\lambda_j - \lambda_k)}e_{jk}& \textrm{if } j>k
\end{array} \right. \end{equation*}
for some $\lambda_1, \ldots, \lambda_n \in \R$ with $\sum_{k=1}^n
\lambda_k = 0$.

By Proposition \ref{hypqc}, it suffices to find a hyper maximal
$q$-corner from $\phi$ to $\imath_\C$.  For this, define $\gamma:
M_{n \times 1}(\C) \rightarrow M_{n \times 1}(\C)$ by

\begin{displaymath} \gamma \left(
\begin{array}{c}
b_1 \\ b_2 \\ \vdots \\ b_n \end{array} \right) = \left(
\begin{array}{c}
\frac{1}{1+i\lambda_1 } b_1 \\
\frac{1}{1+i\lambda_2} b_2 \\
\vdots \\
\frac{1}{1+i\lambda_n} b_n \end{array} \right).
\end{displaymath}
Now define $\Upsilon: M_{n+1}(\C) \rightarrow M_{n+1}(\C)$ by

\begin{displaymath} \Upsilon \left(
\begin{array}{ccc}
A_{n \times n} & B_{n \times 1} \\
C_{1 \times n} & a \end{array} \right) = \left(
\begin{array}{ccc}
\phi(A_{n \times n}) & \gamma(B_{n \times 1}) \\
\gamma^*(C_{1 \times n}) & a
\end{array} \right).
\end{displaymath}
Letting $\lambda_{n+1}=0$, we observe that $\Upsilon$ is the Schur
map satisfying

\begin{equation*} \Upsilon(a_{jk}e_{jk}) = \left\{
\begin{array}{cc}
\frac{a_{jk}}{1+i(\lambda_j - \lambda_k)}e_{jk} & \textrm{if } j<k \\
a_{jk}e_{jk} & \textrm{if } j=k \\
\frac{a_{jk}}{1-i(\lambda_j - \lambda_k)}e_{jk}& \textrm{if } j>k
\end{array} \right. \end{equation*}  for all $j,k=1, \ldots, n+1$ and
$A=(a_{jk}) \in M_n(\C)$.  Since $\sum_{i=1}^{n+1} \lambda_k =
\sum_{i=1}^{n} \lambda_k = 0$, it follows from Theorem \ref{phiu}
that $\Upsilon$ is $q$-positive (in fact, $q$-pure), hence $\gamma$
is a $q$-corner from $\phi$ to $\imath_\C$. Now suppose that
$\Upsilon \geq_q \Upsilon' \geq_q 0$ for some $\Upsilon'$ of the
form

\begin{displaymath} \Upsilon' \left(
\begin{array}{ccc}
A_{n \times n} & B_{n \times 1} \\
C_{1 \times n} & a \end{array} \right) = \left(
\begin{array}{ccc}
\phi'(A_{n \times n}) & \gamma(B_{n \times 1}) \\
\gamma^*(C_{1 \times n}) & \imath'(a)
\end{array} \right).
\end{displaymath}
Since $\Upsilon$ is $q$-pure and $\Upsilon'$ is not the zero map, we
know that $\Upsilon' = \Upsilon^{(t)}$ for some $t \geq 0$, and a
quick calculation gives us

\begin{displaymath} \Upsilon' \left(
\begin{array}{ccc}
A_{n \times n} & B_{n \times 1} \\
C_{1 \times n} & a \end{array} \right) = \left( \begin{array}{cc}
\phi^{(t)}(A_{n \times n}) & \gamma^{(t)}(B_{n \times 1}) \\
(\gamma^*)^{(t)}(C_{1 \times n}) & \frac{1}{1+t}(a)
\end{array} \right).
\end{displaymath}
By inspecting the two formulas for $\Upsilon'$ we see
$\gamma=\gamma^{(t)}$.  But $\gamma^{(t)}$ has the form
\begin{displaymath} \gamma^{(t)} \left(
\begin{array}{c}
b_1 \\ b_2 \\ \vdots \\ b_n \end{array} \right) = \left(
\begin{array}{c}
\frac{1}{1+t+i\lambda_1 } b_1 \\
\frac{1}{1+t+i\lambda_2} b_2 \\
\vdots \\
\frac{1}{1+t+i\lambda_n} b_n \end{array} \right),
\end{displaymath}  hence $t=0$. Therefore, $\Upsilon'=\Upsilon$,
and we conclude the $q$-corner $\gamma$ is hyper maximal. \qed
\end{pf}

In conclusion, we approach the broader question of simply
finding all unital $q$-pure maps $\phi: M_n(\C) \rightarrow M_n(\C)$, as
they provide us with the simplest way to construct and compare
$E_0$-semigroups
through boundary weight doubles.  We believe that all $q$-pure maps
are invertible or have rank one.  For $n=2$, we find in \cite{Me2}
that this conjecture holds: There is no unital $q$-pure map $\phi:
M_2(\C) \rightarrow M_2(\C)$ of rank 2, and there is no unital
$q$-positive map $\phi: M_2(\C) \rightarrow M_2(\C)$ of rank 3. 
It seems that for $n=3$, the key to classifying unital $q$-pure
maps is through investigation of the limits $L_\phi = \lim_{t
\rightarrow \infty} t \phi(I + t \phi)^{-1}$, though the situation
becomes very complicated if $n>3$.

\begin{center} {\bf Acknowledgments} \end{center}
The author very gratefully thanks his thesis advisor, Robert Powers,
for his boundless enthusiasm, constant encouragement, and guidance
in research. His help in the author's thesis work has been
indispensable.  The author would also like to thank Geoff Price for 
proofreading an earlier draft of the paper and making suggestions.

\newpage \section{Notes on Version 2}
The text of this version should be identical to that of the version appearing in the Journal
of Functional Analysis.
\\ \\There is an update for reference \cite{Me2}:
\newline C. Jankowski, \emph{Unital $q$-positive maps
on $M_2(\C)$ and a related $E_0$-semigroup result}, \newline arXiv:1005.4404v1.
\\ \\The author apologizes that some typos remain:
\\ \\In Theorem 2.2, 
the equation of minimality should read ``$\alpha_{t_1}^d(
W A_1 W^*) \cdots \alpha_{t_n}^d(WA_nW^*)Wf$," and we should have ``$f \in H_1$" rather than
``$f \in K$."
\\ \\In the last equation of the proof of Lemma 3.5, the $S_i'$ and $T_i'$ terms should be $\tilde{S}_i$
and $\tilde{T}_i$ terms.  In this same equation, the numbers $a_{ij}$ should be $r_{ij}$.
\\ \\The last line of page 23 should read ``$\ell_t(\rho_{12})(A)= \ell(\rho_{12})\Big(W_t[W_t^*AX_t]X_t^*\Big)=0.$''
\end{document}